\documentclass[leqno]{article}
\usepackage{amssymb,latexsym,amsmath}
\newtheorem{itheorem}{Theorem}
\newtheorem{theorem}{Theorem}[section]

\newtheorem{prop}[theorem]{Proposition}
\newtheorem{corol}[theorem]{Corollary}
\newtheorem{defi}[theorem]{Definition}
\newtheorem{lemma}[theorem]{Lemma}

\def\noproof{{\unskip\nobreak\hfill\penalty50\hskip2em\hbox{}%
     \nobreak\hfill$\square$\parfillskip=0pt%
     \finalhyphendemerits=0\par}}
\def\enddemo{\ifmmode\eqno\square\else\noproof\vskip 12pt plus 3pt minus 9pt
\fi}
\def\diagram{\renewcommand\arraystretch{1.5} $$ \begin{array}}
\def\enddiagram{\end{array} $$ \renewcommand\arraystretch{1}}

\def\rmnewname#1{\expandafter\gdef\csname#1\endcsname{{\mathop{\rm
#1}\nolimits}}}
\def\itnewname#1{{\expandafter\gdef\csname#1\endcsname{{\mathop{\it
#1}\nolimits}}}}

 \itnewname{ab}
 \itnewname{can}
 \rmnewname{cd}
 \itnewname{cl}
 \rmnewname{codim}
 \rmnewname{coker}
 \rmnewname{cone}
 \itnewname{cor}
 \rmnewname{CH}
 \rmnewname{div}
 \rmnewname{Div}
 \itnewname{et}
 \rmnewname{Et}
 \rmnewname{Ext}
 \itnewname{Frob}
 \rmnewname{Hom}
 \itnewname{id}
 \rmnewname{im}
 \rmnewname{Ker}
 \rmnewname{length}
 \rmnewname{Mor}
 \rmnewname{NK}
 \rmnewname{Pic}
 \rmnewname{Quot}
 \itnewname{rec}
 \rmnewname{red}
 \rmnewname{reg}
 \itnewname{res}
 \rmnewname{Sch}
 \itnewname{sh}
 \rmnewname{SK}
 \rmnewname{Sm}
 \rmnewname{SmCor}
 \itnewname{sp}
 \rmnewname{Spec}
 \rmnewname{supp}
 \itnewname{tame}
 \rmnewname{Tame}
 \rmnewname{top}
 \rmnewname{Tor}
 \rmnewname{Tr}
 \itnewname{Zar}
 \def\C{{\mathbb C}}
 
 \def\G{{\mathbb G}}

 \def\M{{\mathcal M}}
 \def\N{{\mathbb N}}
 \def\O{{\mathcal O}}

 \def\p{{\mathfrak p}}

 \def\Q{{\mathbb Q}}
  \def\P{{\mathcal P}}
  
   \def\R{{\mathbb R}}
  \def\X{{\mathfrak X}}
  
  \def\D{{\mathfrak D}}
 \def\Z{{\mathbb Z}}

 \def\noi{\noindent}
 
 \def\ms{\medskip}
 \def\ds{\displaystyle}

 \def\sst{\scriptstyle}
 \def\ssz{\scriptsize}

 \def\Cal#1{{\cal #1}}

 \def\lang{\longrightarrow}
 \def\longeq{\hbox{$=\!=$}}
 \def\mapr#1{\mathrel{\mbox{$\stackrel{#1}{\longrightarrow}$}}}
 
 \def\mapd#1{\Big\downarrow\rlap{$\vcenter{\hbox{$\scriptstyle{{#1}}$}}$}}
 
 \def\liso{\mathrel{\hbox{$\longrightarrow$} \kern-14pt\lower-4pt%
 \hbox{$\scriptstyle\sim$}\kern7pt}}
 \def\r@iso{\mathrel{\lower2pt\hbox{$\scriptstyle\sim$}%
 \kern-8pt\hbox{$\rightarrow$}}}
 \def\riso#1{\mathrel{\stackrel{\!\!#1}{\r@iso}}}
 \def\lr@iso{\mathrel{\kern6pt\lower2pt\hbox{$\scriptstyle\sim$}%
  \kern-12pt\hbox{$\longrightarrow$}}}
 \def\lriso#1{\mathrel{\mathop{\lr@iso}\limits^{#1}}}
 \def\eqd{\Big\|}
 
 \def\mapu#1{\Big\uparrow\rlap {$\vcenter{\hbox{$\scriptstyle{{#1}}$}}$}}
 \def\mapd#1{\Big\downarrow\rlap {$\vcenter{\hbox{$\scriptstyle{{#1}}$}}$}}
 \def\surjr#1{\stackrel{#1}{\hbox{$\relbar \! \! \! \twoheadrightarrow$}}}

\def\surjd#1{\lower4pt\hbox{$\downarrow$}\kern-.57em\Big\downarrow\rlap%
  {$\vcenter{\hbox{$\scriptstyle{{#1}}$}}$}}
 \def\surju#1{\lower-4pt\hbox{$\uparrow$}\kern-7pt\Big\uparrow\rlap%
  {$\vcenter{\hbox{$\scriptstyle{{#1}}$}}$}}

 \def\injd#1{
 \setlength{\unitlength}{0.1pt}
 \begin{picture}(40,120)
 \put(20,102){\vector(0,-1){158}} \put(13.5,100){\mbox{$\scriptscriptstyle
 \cap $}} \put(47,10){\mbox{$\scriptstyle #1 $}}
 \end{picture}}

 \def\dirlim{\mathop{\lim\limits_{{\longrightarrow}}\,}}
 \def\projlim{\mathop{\lim\limits_{{\longleftarrow}}\,}}

\def\tcup{\mathop{\mathord{\cup}\mkern-8.5mu^{\hbox{.}}\mkern 9mu}}

\def\Tcup{\textstyle \tcup}

 \def\hang{\hangindent\Itemindent}
 \def\textindent#1{\hskip\Itemindent\llap{\hbox{\rm #1}\enspace}\ignorespaces}
 \def\Item{\par\noindent\hang\textindent}
 
 \newdimen\Itemindent \Itemindent=.9cm

  \def\m{{\mathfrak m}}
 \def\Zt{{\mathcal Z}}
 \def\Ce{{\mathcal C}}

\author{Alexander Schmidt}
\date{April 26, 2002}
\title{\LARGE \bf \boldmath Relative $K$-groups and class field theory for arithmetic surfaces}
\begin{document}

\maketitle
\section{Introduction} Let $X$ be a regular connected scheme, flat and proper over
$\Spec(\Z)$ whose generic fibre $X_\Q=X \otimes_{\Z}\Q$ is projective over $\Q$. Sending the
class $[P]$ of a closed point $P$ on $X$ to its Frobenius automorphism $\Frob_P$, it is
well-known that one obtains a reciprocity homomorphism
\[
\rec_{X}: \CH_0(X) \lang \tilde \pi_1^\et(X)^\ab.
\]
Here $\CH_0(X)$ denotes the Chow group of zero-cycles on $X$ modulo rational equivalence and
$\tilde \pi_1^\et(X)^\ab$ is the modified abelianized \'{e}tale fundamental group, which
classifies finite abelian \'{e}tale coverings of $X$ in which every real point splits completely
(if $X$ is defined over a totally imaginary number field, then this is just the usual
abelianized \'{e}tale fundamental group). By Lang \cite{La}, the homomorphism $\rec_{X}$ has a
dense image and (after S. Bloch \cite{Bl} solved special cases) K. Kato and S. Saito
\cite{K-S1}, \cite{K-S2}, proved the following theorem which is usually cited as
``unramified class field theory".

\medskip\noindent
{\bf Theorem}  $\rec_{X}$ is an isomorphism of finite abelian groups.

\medskip\noindent
Now assume that $Y$ is the support of a divisor on $X$ and put $U=X-Y$. We consider the
abelianized modified tame fundamental group $\tilde \pi_1^t(X,Y)^\ab$ which classifies finite
abelian \'{e}tale coverings of $U$ which are tamely ramified along $Y$ (see
section~\ref{tamesect} below) and in which every real point splits completely. In \cite{S1},
the author showed the following

\begin{itheorem} \label{finite}
The group $\tilde \pi_1^t(X,Y)^\ab$ is finite and depends only on the scheme $U=X-Y$ (and not
on the choice of\/ the compactifying scheme $X$).
\end{itheorem}

Recall that the group $\CH_0(X)$ can be identified with the spectral term $E_2^{d,-d}$
($d=\dim X$) of the Quillen spectral sequence
\[
E_1^{pq}= \bigoplus_{x\in X^p} K_{-p-q}(k(x))\Rightarrow K_{-p-q}(X),
\]
in which $K_*(-)$ are Quillen's $K$-groups \cite{Qu}. There exists an analogous spectral
sequence for relative $K$-theory and we define the relative Chow group of zero-cycles
$\CH_0(X,Y)$ as the $E_2^{d,-d}$-term of that spectral sequence (see section~\ref{relk}
below).

\medskip
The objective of this paper is to give a proof of the

\begin{itheorem} \label{main} If\/ $X$ is an arithmetic surface (i.e.\ $\dim_{Krull} X=2$), then
there exists a natural isomorphism of finite abelian groups
\[
\rec_{X,Y}: \CH_0(X,Y) \lang \tilde \pi_1^t(X,Y)^\ab\;.
\]
\end{itheorem}

\noindent {\bf Remarks.}
1. It is natural to conjecture, that the statement of the above theorem holds in any
dimension. This is easily shown in dimension~$1$ and we think that in principle the same
methods as in dimension~$2$ can also be applied to any dimension $>2$. However the technical
effort will increase rapidly.

\smallskip\noindent
2. By theorem~\ref{finite}, the right hand side of the isomorphism $\rec_{X, Y}$ only depends
on $U=X-Y$ and not on the choice of $X$. Hence the same is also true for the left hand side,
and it is desirable to give an intrinsic definition of it only in terms of $U$. The author
conjectures that $\CH_0(X,Y)$ is naturally isomorphic to the 0th singular homology group
$h_0(U)$, which is considered in \cite{S2}. This is true in the analogous situation for
varieties over finite fields and in that situation the corresponding statement to that of
theorem~\ref{main} is known in any dimension, see \cite{S-S}. In the mixed characteristic
case, there exists a natural surjection $h_0(U) \twoheadrightarrow \CH_0(X,Y)$ and a small
but unfortunately highly resistant moving lemma {\em would} show that this map is an
isomorphism.

\smallskip\noindent
3. In \cite{K-S3}, Kato and Saito give an interpretation of $\pi_1(U)^\ab$ in terms of a
higher id\`{e}le class group which is defined as a Nisnevi\v c cohomology group of a relative
Milnor $K$-sheaf. One can interpret theorem~\ref{main} as a duality result on relative
$K$-groups.

\section{Results from the work of Kato and Saito}\label{classfieldsect}

For reference, to fix notation and to deduce some  corollaries,  we recall several results
from the work of K. Kato and S. Saito. Let $X$ be a connected noetherian scheme. For $j \in
\N$, let
\[
X_j=\left\{x \in X\,|\, \dim \left(\overline{\{x \}}\right)=j \right\}.
\]
For an integer $i\geq 0$, the localization theory in Quillen $K$-theory on $X$ gives rise to
a homomorphism
\begin{equation}\label{divmap}
 \delta_i: \bigoplus_{y \in X_1} K_{i+1}(k(y)) \lang \bigoplus_{x \in X_0} K_i(k(x)).
\end{equation}
We define $\SK_i(X)$ to be the cokernel of this map $\delta_i$. Assume that $X$ is a variety
over a field $k$ and that $k$ (and hence also every finite extension of $k$) carries a
natural topology. Then we give $\SK_0(X)$ the discrete topology and we endow $\SK_1(X)$ with
the finest topology such that for every $x\in X_0$ the natural homomorphism
\[
k(x)^\times \lang \SK_1(X)
\]
is continuous.

\bigskip\noindent
\underline{1. The reciprocity homomorphism for an arithmetic scheme.}

\medskip Let $X$ be a connected scheme of finite type over $\Spec(\Z)$. For each closed point $x
\in X_0$, $k(x)$ is a finite field, so there is an isomorphism
\[
\hat \Z \liso \pi_1^\ab(x)=Gal(k(x)^{sep}|k(k)),
\]
which sends $1\in \hat \Z$ to the Frobenius $f_x$ over $k(x)$. We define the Frobenius
element $F_x$ of $x$ in $\pi_1^\ab(X)$ to be the image of $f_x$ under the natural
homomorphism $\pi_1^\ab(x) \to \pi_1^\ab(X)$. The assignment $1 \mapsto F_x$, $x \in X_0$,
defines a homomorphism

\begin{equation} \label{recmap}
\bigoplus_{x \in X_0} \Z \lang \pi_1^\ab (X).
\end{equation}

Then by \cite{La}, we have the
\begin{lemma}
If $X$ is irreducible and the reduced subscheme $X_{red}$ of $X$ is normal, then the map
\eqref{recmap} has a dense image.
\end{lemma}

Now we consider the abelianized modified fundamental group $\tilde\pi_1^\ab(X)$ which is the
unique quotient of $\pi_1^\ab(X)$ classifying finite abelian \'{e}tale coverings of $X$ in which
every real-valued point of $X$ splits completely. Then we have the

\begin{lemma}{\rm (\cite{Sa}, lemma 2.4)}
If $X$ is proper over $\Spec(\Z)$, then the composite of the map \eqref{recmap} and the
natural surjection $\pi_1^\ab (X) \to \tilde \pi_1^\ab (X)$ annihilates the image of
\[
\delta: \bigoplus_{x \in X_1} k(x)^\times \lang \bigoplus_{x \in X_0} \Z,
\]
which is the map \eqref{divmap} for $i=0$.
\end{lemma}

Consequently, we obtain a natural map
\[
\SK_0(X) \lang \tilde \pi_1^\ab(X),
\]
which is called the reciprocity map for $X$.

\bigskip\noindent
\underline{2. Class field theory for varieties over $\R$ and $\C$.}

\medskip
We follow \cite{Sa}, \S4. Let $k^*$ be $\R$ or $\C$ and let $X^*$ be a connected, proper and
smooth scheme over $k$. For $x \in X^*_0$, $k(x)\cong \R$ or $\C$ and we have a canonical
surjection
\[
k(x)^\times \lang \pi_1^\ab(x)=Gal({k(x)}^{sep}|k(x)).
\]
The maps $\pi_1^\ab(x) \to \pi_1^\ab(X*)$ for $x \in X^*_0$ induce a well-defined continuous
homomorphism
\begin{equation}\label{eq1}
\tau: \SK_1(X^*) \lang \pi_1^\ab(X^*),
\end{equation}
which is the zero map if $k^*=\C$ or $X^*(\R)=\varnothing$.  Let $k\subset k^*$ be a
subfield and suppose that
\Itemindent=.8cm

 \Item{-} $k$ is dense in $k^*$ for the usual topology of $k^*$.
 \Item{-} If $k^*=\C$, $k$ is algebraically closed.
 If $k^*=\R$, the algebraic closure $\bar k$
of $k$ is a quadratic extension of $k$.
 \Item{-} There is a proper, smooth scheme $X$ over $k$ such that
\[
X \otimes_k k^* \cong X^*.
\]

By the proper base change theorem, the natural map $\pi_1^\ab(X^*)\to \pi_1^\ab(X)$ is an
isomorphism and the image of the natural homomorphism
\begin{equation}\label{eq2}
\SK_1(X) \lang \pi_1^\ab(X)
\end{equation}
coincides with that of the map (\ref{eq1}). Slightly more general, let $V \subset X$ be a
nonempty open subscheme and $V^*=V\otimes_k k^*$. Using the smooth base change theorem
instead (we are in characteristic zero), we have again the isomorphism $\pi_1^\ab(V^*) \liso
\pi_1^\ab(V)$ and, sending a point $x \in V(\R)=V^*(\R)$ to the image of the unique
nontrivial element of $\pi_1^\ab(x)\cong Gal(\C|\R)$ in $\pi_1^\ab(V)$, we obtain a natural
map
\[
i: V(\R) \lang \pi_1^\ab(V).
\]
The map $i$ is locally constant on $V(\R)$ and has finite image  (see \cite{Sa}, lemma (4.8)
for the case $V=X$, the proof in the general case is the same). If $\chi\in
H^1_\et(V,\Q/\Z)$ corresponds to a cyclic \'{e}tale covering $\tilde V \to V$, then a point $x
\in V(\R)$ splits completely in $\tilde V$ if and only if $\chi$ is trivial on $i(x)$. In
particular, the subset of points in $V(\R)$ which split completely in $\tilde V$ is (norm)
open and closed in $V(\R)$.

\bigskip\noindent
\underline{3. Class field theory of schemes over henselian discrete valuation fields.}

\medskip
We follow \cite{Sa}, \S3. Let $\O_k$ be a henselian discrete valuation ring with finite
residue field $F$ and quotient field $k$. Let $X$ be a proper smooth scheme over $k$. For $x
\in X_0$, $k(x)$ is a finite extension of $k$, so that the local class field theory for
$k(x)$ gives us a canonical homomorphism
\[
\rho_x : k(x)^\times \lang \pi_1^\ab(x) \lang \pi_1^\ab (X).
\]
The sum of these $\rho_x$ is trivial on the diagonal image of $\bigoplus_{x \in
X_1}K_2(k(X))$ and we obtain a canonical morphism
\begin{equation}\label{eq3}
 \tau: \SK_1(X) \lang \pi_1^\ab(X).
\end{equation}
Giving each $k(x)$, $x\in X_0$ the usual adic topology, we obtain a natural topology on
$\SK_1(X)$ (see the beginning of this section). The map  $\tau$ is continuous with respect to
this topology and the natural profinite topology on $\pi_1^\ab(X)$. For every natural number
$n$ which is prime to char$(k)$, the subgroup $n\SK_1(X)$ is open in $\SK_1(X)$. In
particular, if char$(k)$=0, then every homomorphism $\SK_1(X) \to \Q/\Z$ with finite image is
automatically continuous.

Now assume that $X$ has a proper regular model $\X$ over $\O_k$. By \cite{Sa}, lemma 3.1, we
have a natural homomorphism $\delta: \SK_1(K) \to \SK_0(Y)$ and a commutative diagram
 \diagram{ccc}
 \SK_1(X)&\lang & \pi_1^\ab(X)\\
 \mapd{\delta}&&\mapd{sp}\\
 \SK_0(Y) &\mapr{\delta} &\pi_1^\ab (Y),
 \enddiagram
where $sp$ is the specialization homomorphism (cf.\ \cite{SGA1}, X) and $\phi$ is the
reciprocity map for $Y$ as defined in paragraph 1 above. Important for us is the

\begin{prop} {\rm (\cite{Sa}, prop. 3.3)}
Let $\chi \in H^1_\et(X,\Q/\Z)$ and let $\tilde \chi: \SK_1(X) \to \Q/\Z$ be the induced
homomorphism via \eqref{eq3}. Then $\chi$ comes from the subgroup $H^1_\et(Y,\Q/\Z)\cong
H^1_\et(\X,\Q/\Z) \subset H^1_\et(X,\Q/\Z)$ if and only if $\tilde \chi$ factors though the
map $\delta$.
\end{prop}

The crucial step in the proof of the last proposition is the following lemma, which we will
need in the following.

\begin{lemma}{\rm (\cite{Sa}, lemma 3.15)} \label{saitolemma}
Let $A$ be a henselian regular local ring of dimension
$\geq 2$ with perfect residue field $F$ and quotient field $K$. If\/ $\hbox{\rm char}(K)=0$
and $\hbox{\rm char}(F)=p>0$ assume that there exists exactly one prime ideal $\p$ of height
$1$ which divides $p$ and let $T\in A$ with $\p=(T)$. In the remaining cases let $T\neq 0$ be
any non-unit of $A$. Put $U=\Spec(A[1/T])$. Then, if $\chi \in H^1_\et(U,\Q/\Z)$ induces an
unramified character $\chi_u \in H_\et^1(u,\Q/\Z)$ for each $u \in U_0$, $\chi$ comes from
$H^1_\et(\Spec(A),\Q/\Z)$.
\end{lemma}

\bigskip\noindent
\underline{4. Unramified class field theory of arithmetic schemes.}

\medskip\noindent
Now let $\X$ be a connected regular, proper and flat scheme over $\Spec(\Z)$, and suppose
that $X= \X \otimes_\Z \Q$ is projective over $\Spec(\Q)$.  Let $k$ be the algebraic closure
of $\Q$ in the function field of $X$ and put $S=\Spec(\Cal O_k)$. The structural morphism $\X
\to \Spec(\Z)$ factors through $S$ and $X$ is geometrically irreducible as a variety over
$k$. Let $S_f$ be the set of closed points of $S$ and let $S_\infty$ be the set of
archimedean places of the number field $k$. For $v \in S_f$ let $Y_v = \X \otimes_S v$ be
the special fibre of $\X$ over $v$. For $v \in S_\infty \cup S_f$ let $k_v$ be the algebraic
closure of $k$ in the completion of $k$ at $v$ and $X_v=X \times_k k_v$.

\ms
For $v \in S_f \cup S_\infty$ we endow the fields $k_v$ with the restriction of the
natural topology of the completion of $k$ at $v$ to $k_v$. For $v \in S_f$ there exists a
natural, surjective and continuous (\cite{Sa}, 3.11) homomorphism
\[
\SK_1(X_v) \lang \SK_0(Y_v)=\CH_0(Y_v).
\]
For an open subscheme $U \subset S$ we consider the topological group
\[
I(\X/U)= \left(\prod_{v \in P_U} \SK_1(X_v) \right) \times \left(\bigoplus_{v\in U_0}
\CH_0(Y_v) \right) ,
\]
where $P_U$ denote the set of places (including the archimedean ones) which are not in $U_0$.
There exists a natural homomorphism (\cite{Sa}, 5.3)
\[
I(\X/U) \lang \pi_1(\X_U)^{\ab}
\]
which annihilates the image of $\SK_1(X) \lang I(\X/U)$, inducing a homomorphism (\cite{Sa},
5.5)
\[
\tau: C(\X/U):=\coker\big(\SK_1(X) \to I(\X)\big) \lang \pi_1(\X)^\ab.
\]
By \cite{Sa}, 5.6., every subgroup of finite index in $C(\X/U)$ is open and, for every
positive integer $n$, $\tau$ induces an isomorphism of finite abelian groups
\[
C(\X/U)/n \mapr{\sim} \pi_1(\X_U)^\ab /n.
\]

\bigskip\noindent
\underline{5. Class field theory for two-dimensional henselian fields.}

\medskip Let $F$ be a discrete valuation field with residue field $k$. Let us recall the
filtration on $K_r(F)$ for $r=1,2$. Let $\O_F$ be the valuation ring in $F$ and let $\m
\subset \O_F$ be its maximal ideal, thus $\O_F/\m \cong k$. Then one puts
$U^0(F^\times)=\O_F^\times$ and for $i \geq 1$, $U^i(F^\times)=\ker(\O_F^\times \to
(\O_F/\m^i)^\times)$. As is well known, we have isomorphisms
\[
F^\times/U^0F^\times \cong \Z\;, \qquad U^0(F^\times)/U^1(F^\times)\cong k^\times.
\]
For $i\geq 1$ the group $U^iK_2(F)$ is the subgroup generated by symbols $(u,x)$ with $u\in
U^i(F^\times)$ and $x \in F^\times$. The group $U^0K_2(F)$ is the kernel of the tame symbol
$K_2(F) \to k^\times$. We have inclusions
\[
K_2(F) \supseteq U^0K_2(F) \supseteq U^1K_2(F) \supseteq \cdots
\]
and natural isomorphisms (\cite{B-T}, prop.4.3,4.5)
\[
K_2(F)/U^0K_2(F)\cong k^\times\; , \qquad U^0K_2(F)/U^1K_2(F) \cong K_2(k).
\]

Now let $K$ be a henselian discrete valuation field whose residue field $F$ is a henselian
discrete valuation field with finite residue field. Such fields $K$ are called
two-dimensional henselian fields. The following results are contained in \cite{Ka},
\cite{K-S2} for complete fields, but the results carry over to the henselian case.

We consider the group $K_2(K)$ and the filtration $U^iK_2(K)$ which we get by considering
$K$ just as a discrete valuation field. For $i\geq 0$ there exist compatible topologies on
$K_2(K)/U^iK_2(K)$ (see \cite{K-S2}, \S2). These topologies might depend on some choices
(cf.\ the discussion in loc.cit.), but at least the topology on $K_2(K)/U^1K_2(K)$ is
independent of these choices. Assume for simplicity that ${\rm char}(K)=0$ (we will only
need this case). Then, setting for $q\geq 1$
\[
H^q(K):=H^q(K,\Q/\Z(q-1)),
\]
there exists a natural pairing
\[
H^1(K) \times K_2(K) \lang H^3(K) \liso \Q/\Z
\]
which induces a homomorphism $H^1(K) \to \Hom(K_2(K),\Q/\Z)$. With respect to this
homomorphism, we have the

\begin{theorem} \begin{description}
\item{\rm (i)} If ${\rm char}(F)=p> 0$, then $H^1(K)$ is isomorphic to the group of
all homomorphisms $\phi: K_2(K)\to \Q/\Z$ such that $\phi(U^iK_2(K))=0$ for some $i$ and
such that the induced homomorphism $K_2(K)/U^iK_2(K) \to \Q/\Z$ is continuous with respect
to the discrete topology on $\Q/\Z$. For each $i \geq 1$, the abelian group
$U^iK_2(K)/U^{i+1}K_2(K)$ is annihilated by $p$.
\item{\rm (ii)} If ${\rm char}(F)= 0$, then $H^1(K)$ is isomorphic to the group of
all homomorphisms $\phi: K_2(K)\to \Q/\Z$ of finite order. The group $U^1K_2(K)$ is
divisible and annihilated by any homomorphism $K_2(K)\to \Q/\Z$ of finite order. For $i\geq
1$, any homomorphism $K_2(K)/U^iK_2(K)\to \Q/\Z$ of finite order is continuous with respect
to the discrete topology on $\Q/\Z$.
\end{description}

\medskip\noindent In both cases, a homomorphism $\phi: K_2(K) \to \Q/\Z$ corresponding to an
element in $H^1(K)$ is unramified, i.e.\ lies in the subgroup $H^1(F,\Q/\Z)\subset H^1(K)$,
if and only if $\phi(U^0K_2(K))=0$.
\end{theorem}

We conclude the following

\begin{corol}
Let $\chi\in H^1(K,\Q/\Z)$ and let $\tilde\chi: K_2(K) \to \Q/\Z$ be the corresponding
homomorphism. Then the cyclic field extension $K_\chi|K$ corresponding to $\chi$ is tamely
ramified (with respect to the structure of $K$ as a usual discrete valuation field) if and
only if $\tilde \chi(U^1K_2(K))=0$.
\end{corol}

\bigskip\noindent
\underline{6. Class field theory for two-dimensional henselian rings.}

\medskip
We follow \cite{K-S2}, \S2, \S4, where the complete case is considered, but all arguments
carry over to the henselian case. Let $A$ be a two-dimensional regular local henselian ring
with finite residue field $k$ and quotient field $K$. Again, we assume for simplicity that
${\rm char}(K)=0$. Let $P$ be the set of prime ideals of height one in $A$. For $z\in P$ we
denote by $K_z$ the henselization of $K$ at $z$, which is a two-dimensional henselian field.
A modulus $m$ is a family $(m(z))_z$ of nonnegative integers given for each point $z\in P$
such that $m(z)=0$ for all but finitely many $z$. For a modulus $m$, let
\[
I_m(K) = \bigoplus_{z \in P} K_2(K_z)/U^{m(z)}K_2(K_z)
\]
\[
C_m(K)=\coker(K_2(K) \mapr{\rm diag} I_m(K)).
\]
The complex of Bloch-Gersten-Quillen
\[
K_2(K) \to \bigoplus_{z \in P} K_1(k(z)) \to K_0(k) \to 0
\]
and the tame symbols $K_2(K_z) \to K_1(k(z)) $ give a canonical homomorphism $d: C_m(K) \to
K_0(k)=\Z$, which is an isomorphism if $m=(0)$. Let $$C(K)= \projlim_m C_m(K)$$ where $m$
ranges over all moduli. We endow $C_m(K)$ with the discrete topology and $C(K)$ with the
topology as an inverse limit. For $\chi\in H^1(K)$, its image $\chi_z$ in $H^1(K_z)$, $z\in
P$, induces a homomorphism $K_2(K_z)/U^iK_2(K_z)$ for some $i$ which depends on $z$. Since
$\chi_z$ is unramified for almost all $z$, there exists a modulus $m$ such that $\chi$
defines a homomorphism $I_m(K)\to \Q/\Z$. One can show that this homomorphism factors
through $C_m(K)$ and we obtain a natural homomorphism
\[
H^1(K) \lang \Hom(C(K), \Q/\Z).
\]

\begin{theorem}
The above homomorphism is injective and its image consists of all continuous homomorphisms
with finite image.
\end{theorem}

Now let $z_1,\ldots,z_r \in P$ be given and let $m=m_{z_1,\ldots,z_r}$ be the modulus with
$m(z_1)=\cdots =m(z_r)=1$ and $m(z)=0$ for $z\notin \{z_1,\ldots,z_r\}$. We obtain the
following

\begin{corol} \label{localtamecrit}
Let $z_1,\ldots,z_r \in P$ be given and let $m=m_{z_1,\ldots,z_r}$ be the modulus with
$m(z_1)=\cdots =m(z_r)=1$ and $m(z)=0$ for $z\notin \{z_1,\ldots,z_r\}$. Then the elements\/
$\chi \in H^1(K)$ for which the finite cyclic extension $K_\chi|K$ is unramified at all
$z\in P - \{z_1,\ldots,z_r\}$ and tamely ramified at (the discrete valuations associated
with) $z_1,\ldots,z_r$ correspond bijectively to the set of homomorphisms
\[
C_m(K) \lang \Q/\Z
\]
with finite images.
\end{corol}

\section{Relative \boldmath $K$-groups} \label{relk}

Let us collect some facts on relative $K$-groups. A more detailed description and proofs of
several facts mentioned below are found in \cite{Le}, \S1.

\medskip
Following Quillen \cite{Qu}, the $K$-groups of an exact category $\Ce$ are defined by
\[
K_i\Ce=\pi_{i+1} (K\Ce,*), \quad K\Ce= BQ\Ce.
\]
Here $B$ means geometric realization, $Q$ is Quillen's construction \cite{Qu} and
$\pi_i(-,*)$ is the $i$-th homotopy group of a pointed topological space. For a scheme $X$
let $\M_X$ be the category of coherent $\O_X$-sheaves and $\P_X$ the full subcategory of
locally free $\O_X$-sheaves. As usual, we set $K_iX= K_iQ\P_X$, and $G_iX=K_iQ\M_X$. The
exact inclusion $\P_X \to \M_X$ induces homomorphisms
\[
K_iX \lang G_iX,
\]
which are isomorphisms when $X$ is regular. If $Y$ is a closed subscheme of $X$ then $j_Y^*:
\P_X \to \P_Y$ is exact and one defines the relative $K$-groups $K_*(X,Y)$ as the homotopy
groups of the homotopy fibre $K(X,Y)$ of the induced map $j_Y^*: KX \to KY$. In particular,
we obtain a long exact sequence
\[
\cdots \to K_{i+1} Y  \to K_i(X,Y) \to K_iX \to K_iY \to \cdots.
\]
Let $\M_{(X,Y)}$ be the full subcategory in $\M_X$ consisting of coherent $\O_X$-sheaves $F$
with $\Tor^{\O_X}_i(F,\O_Y)=0$ for $i>0$. Then the natural functor $j_Y^*: \M_{(X,Y)} \to
\M_Y$ is exact and one defines the relative $G$-groups $G_*(X,Y)$ as the homotopy groups of
the homotopy fibre of $j_Y^*: K\M_{(X,Y)}\to K\M_Y$. If $\O_Y$ has finite
$\O_X$-$\Tor$-dimension (i.e.\ if $X$ is regular or if $Y$ is locally principal), then the
inclusion $\M_{(X,Y)} \to \M_X$ induces an isomorphism on $K$-groups. In this case we obtain
a commutative ladder
\[
\renewcommand\arraystretch{1.5}
\begin{array}{ccccccccc}
 \cdots & \to & K_i(X,Y)&\to & K_i(X) & \to & K_i(Y) & \to & \cdots\\
 &&\mapd{}&&\mapd{}&&\mapd{}\\
 \cdots & \to & G_i(X,Y)&\to & G_i(X) & \to & G_i(Y) & \to & \cdots
\end{array}
\renewcommand\arraystretch{1}
\]
and thus, if $X$ and $Y$ are regular, we obtain isomorphisms
\[
K_i(X,Y)\liso G_i(X,Y).
\]
Let $f:X' \to X$ be a morphism of schemes and let $Y'$ be a closed subscheme in $X'$
contained in $f^{-1}(Y)$. Then we obtain a natural homomorphism
\[
f^*: K_*(X,Y) \lang K_*(X',Y')
\]
and, if $f$ is flat, a similar map $G_*(X,Y)\to G_*(X',Y')$ (see \cite{Co}, \cite{Le} for
the technical details).

\bigskip
From now on we assume that $X$ is a connected two-dimensional regular scheme such that the
residue fields of all closed points are finite. Let $D$ be a divisor on $X$ and $Y=\supp(D)$.
$Y$ is a locally principal subscheme of $X$ defined locally by a non-zero divisor.  Thus, if
$Z$ is a  closed subscheme in $X$, then $\O_Z$ is in $\M_{(X,Y)}$ if and only if $Z$
intersects $Y$ properly. For $F \in \M_{(X,Y)}$ and $\supp(F)=Z$, we have $\O_Z \in
\M_{(X,Y)}$.

\medskip
Let $\M^i_{(X,Y)}$ and $\M_Y^i$ be the subcategories of coherent sheaves $F$ with\linebreak
$\codim_X\supp(F)\geq i$ in $X$ and $Y$, respectively. Let $G(X^i,Y^i)$ be the homotopy
fibre of $BQj_Y*: BQ\M_{(X,Y)}^i \to BQ\M_Y^i$. Then the map
\[
\dirlim_{\sst \renewcommand{\arraystretch}{0.8}\begin{array}{c} \sst Z \subset X\\
\scriptstyle \O_Z\in \M_{X,Y}\\ \sst \codim_XZ \geq i
\end{array}
\renewcommand{\arraystretch}{1}}
G(Z, Z\cap Y) \lang G(X^i,Y^i)
\]
is a homotopy equivalence. For $i <k$, let $\M_{(X,Y)}^{i/k}$ be the direct limit
\[
\M_{(X,Y)}^{i/k}=
\dirlim_{\sst \renewcommand{\arraystretch}{0.8}\begin{array}{c} \sst Z \subset X\\
\scriptstyle \O_Z\in \M_{X,Y}\\ \sst \codim_XZ \geq k
\end{array}
\renewcommand{\arraystretch}{1}}
\M_{(X-Z,Y-Z)}^i
\]
and let $\M_{Y}^{i/k}$ be the direct limit
\[
\M_Y^{i/k}= \dirlim_{\sst \renewcommand{\arraystretch}{0.8}\begin{array}{c} \sst Z \subset Y\\
\sst \codim_X Y \geq k
\end{array}
\renewcommand{\arraystretch}{1}}
\M_{(Y-Z)}^i.
\]
Let $G(X^{i/k},Y^{i/k})$ be the homotopy fibre of $BQ\M^{i/k}_{(X,Y)}\to BQ\M^{i/k}_Y$. Then
the localization sequence for relative $G$-theory (\cite{Le}, 1.5) shows that
\[
G(X^k,Y^k) \to G(X^i,Y^i) \to G(X^{i/k},Y^{i/k})
\]
is a homotopy fibre sequence and one obtains a spectral sequence
\[
E_1^{pq}(X,Y) \Longrightarrow G_{-p-q}(X,Y),
\]
in which the $E_1$-terms are given by
\[
E_1^{pq}(X,Y)=\left\{
\begin{array}{ll}
G_{-p-q}(X^{p/p+1},Y^{p/p+1}),& -p-q>0,p\leq \dim X\\
\bar G_0(X^{p/p+1},Y^{p/p+1}), &-p-q=0\\
\; 0 & \hbox{ otherwise.}
\end{array}
\right.
\]
Here $\bar G_0(X^{p/p+1},Y^{p/p+1})= \im \left(G_0(X^p,Y^p) \to
G_0(X^{p/p+1},Y^{p/p+1})\right)$. The filtration on $G_*(X,Y)$ is the ``topological''
filtration:
\[
F^pG_*(X,Y)= \im \left(  G_*(X^p,Y^p) \lang G_*(X,Y)\right).
\]

\begin{defi} We call the group
\[
\CH^p(X,Y)=E_2^{p,-p}(X,Y)
\]
the relative Chow group of codimension $p$ cycles. If\/ $X$ is equidimensional of dimension
$d$, we call the group
\[
\CH_p(X,Y)=E_2^{d-p,-d+p}(X,Y)
\]
the relative Chow group of $p$-cycles.
\end{defi}

\noindent
{\bf Remark:}  If \/ $Y$ is empty, then
\[
\CH^p(X,\varnothing) = \coker\Big(\bigoplus_{x\in X^{p-1}} k(x)^\times \to \bigoplus_{x\in
X^p} \Z\Big)
\]
is the usual group of codimension $p$ cycles on $X$ modulo rational equivalence.

\medskip
Furthermore, we consider the spectral sequence
\[
E_1^{pq}(X|Y) \Longrightarrow G_{-p-q}(X),
\]
in which the $E_1$-terms are given by
\[
E_1^{pq}(X|Y)=\left\{
\begin{array}{ll}
K_{-p-q}(\M^{p/p+1}_{(X,Y)}),& -p-q>0,p\leq \dim X\\
\bar K_0(\M^{p/p+1}_{(X,Y)}), &-p-q=0\\
\; 0 & \hbox{ otherwise,}
\end{array}
\right.
\]
where $\bar K_0(\M^{p/p+1}_{(X,Y)})= \im \left( K_0(\M^{p}_{(X,Y)}) \to
K_0(\M^{p/p+1}_{(X,Y)})\right)$. Finally we have the usual Quillen spectral sequence for $X$
(and the analogous one for $Y$)
\[
E_1^{pq}(X)= \bigoplus_{x\in X^p} K_{-p-q}(k(x))\Rightarrow G_{-p-q}(X).
\]
The inclusion $\M_{(X,Y)}^{p/p+1} \to \M_X^{p/p+1}$ gives a map on $E_1$-terms which is
compatible with differentials. If the maps $G_0(X^p,Y^p)\to G_0(X^{p/p+1},Y^{p/p+1})$ and
$K_0(\M^{p}_{(X,Y)}) \to K_0(\M^{p/p+1}_{(X,Y)})$ are surjective, then we get a long exact
sequence of $E_1$-terms:
\[
\to E_1^{p,q-1}(Y) \to E_1^{pq}(X,Y) \to E_1^{pq}(X|Y) \to E_1^{pq}(Y) \to
\]
compatible with the differentials. In addition, we consider the  categories \medskip
\Itemindent=2cm
\Item{$\M_{(X|Y)}^i -$}  coherent $\O_X$-sheaves whose support is of codimension $\geq i$ and has
proper intersection with $Y$.

\medskip\noindent
Then we have a natural inclusion $\M_{(X,Y)}^i \to \M_{(X|Y)}^i$.

\begin{lemma}
The natural map $Q\M_{(X,Y)}^i \to Q\M_{(X|Y)}^i$ is a homotopy equivalence for all $i$.
\end{lemma}

\begin{demo}{Proof:}
The statement of obvious for $i = 2$. For $i= 0$ it follows from the regularity of $X$,
since every coherent $\O_X$-sheaf has a resolution by locally free ones. It remains the case
$i=1$. If $F$ is a coherent $\O_X$-sheaf whose support is of codimension $\geq 1$ and has
proper intersection with $Y$, then there exists a closed (not necessarily reduced) subscheme
$Z$ of codimension $1$ in $X$ with $\O_Z \in \M_{(X,Y)}$ such that $F$ is an $\O_Z$-module.
The condition $\O_Z \in \M_{(X,Y)}$ means that $\O_Z$ has no $I_Y$-torsion. Resolving $F$ by
a free $\O_Z$-module, we obtain an exact sequence
\[
0 \to K \to (\O_Z)^n \to F \to 0.
\]
We see that also $K \in \M_{(X,Y)}^1$ and the Quillen's resolution lemma \cite{Qu} implies
the statement for $i=1$.
\end{demo}

For a closed subscheme $Z\subset X$, let $\M_{X}(Z)$ be the category of coherent
$\O_X$-sheaves with support in $Z$. Let $\M_{(X|Y)}^{p/p+1}$ be the direct limit
\[
\M_{(X|Y)}^{p/p+1}=
\dirlim_{\sst \renewcommand{\arraystretch}{0.8}\begin{array}{c} \sst W \subset Z \subset X\\
\scriptstyle \O_W, \O_Z\in \M_{(X,Y)}\\ \sst \codim_XZ \geq p\\
\sst \codim_X W \geq p+1
\end{array}
\renewcommand{\arraystretch}{1}}
\M_{(X-W)}(Z)
\]
The same reasoning as above shows the

\begin{lemma}
The natural map $Q\M_{(X,Y)}^{p/p+1} \to Q\M_{(X|Y)}^{p/p+1}$ is a homotopy equivalence.
\end{lemma}

Let us calculate some terms of the above spectral sequences. This is done in \cite{Le} for
schemes over a field. Using the Gersten conjecture in that range where it is known, we obtain
a proof of the following lemma in exactly the way as in the proof of \cite{Le}, lemma 1.2.

\begin{prop}\label{ausrech}
The map $i_{X-Y}^*: \M_{(X,Y)}^{p/p+1} \to \M_{(X-Y)}^{p/p+1}$ induces an isomorphism
\[
K_0(\M_{(X,Y)}^{p/p+1}) \liso K_0(\M_{(X-Y)}^{p/p+1}) \liso \bigoplus_{x \in (X-Y)^p} \Z \; .
\]
For $q= 1,2$, we obtain short exact sequences

\diagram{ccccccccc}
 0 &\to& K_q(\M_{(X,Y)}^{p/p+1})& \to& K_q(\M_{(X-Y)}^{p/p+1})& \stackrel{\delta}{\to}&
K_{q-1}(\M_{Y}^{p/p+1}) &\to &0\\
&&&&\mapu{\wr}&&\mapu{\wr} &&\\
&&&&\ds\bigoplus_{x \in (X-Y)^p}\!\!\!\!K_q(k(x))&&\ds\bigoplus_{y \in Y^p} K_{q-1}(k(y)).
\enddiagram
\end{prop}

\begin{demo}{Proof:} By dimension reasons, we have an isomorphism
\[
Q\M_{(X|Y)}^{2/3} \to Q\M_{X-Y}^{2/3}.
\]
The localization theorem shows that
\[
Q\M_{Y}^{p/p+1} \to Q\M_{(X|Y)}^{p/p+1} \to Q\M_{X-Y}^{p/p+1}
\]
is a homotopy fibre sequence for $p=0,1$. Let $x \in Y$ be a closed point. Choose a closed
irreducible reduced subscheme $Z$ of codimension $1$ in $X$ which has proper intersection
with $Y$ and such that $x$ is a regular point of $Z$. Let $R$ be the semi-local ring of
$Y\cap Z$ in $Z$ and denote its normalization by $\tilde R$. Then $\tilde R$ is a semi-local
PID with finite residue fields and by Gersten \cite{Ge} the map
\[
K_*(k(x)) \to K_*(\tilde R) \to G_*(R)
\]
is the zero-map. Hence the map $K_*\M_{Y}^{1/2} \to K_*\M_{(X|Y)}^{1/2}$ is zero. The same
reasoning, but now using the result of Dennis-Stein \cite{DS}, shows that the map
\[
K_q(\M_Y^{0/1}) \to K_q(\M_{(X|Y)}^{0/1})
\]
is zero for $q \leq 2$.  Summing up, we obtain the required statements.
\end{demo}

\begin{prop}\label{surj}
The natural maps
\[
K_0(\M_{(X,Y)}^{p}) \to K_0(\M_{(X,Y)}^{p/p+1})\quad \hbox{and} \quad G_0(X^p,Y^p) \to
G_0(X^{p/p+1},Y^{p/p+1})
\]
are surjective. In particular, we get a long exact sequence of $E_1$-terms:
\[
\to E_1^{p,q-1}(Y) \to E_1^{pq}(X,Y) \to E_1^{pq}(X|Y) \to E_1^{pq}(Y) \to
\]
compatible with the differentials.
\end{prop}

\begin{demo}{Proof:} By the last proposition, we have
\[
K_0(\M_{(X,Y)}^{p/p+1}) \liso \bigoplus_{x \in (X-Y)^p} \Z  \cong \bigoplus_{x \in (X-Y)^p}
K_0(k(x))\;.
\]
Similarly,
\[
K_0(\M_Y^{p/p+1}) \liso\bigoplus_{x \in Y^p} \Z \cong \bigoplus_{x \in Y^p} K_0(k(x))\;.
\]
In the commutative ladder with exact rows

\diagram{cccccc}
 \to & G_0(X^p, Y^p)&\to &K_0(\M_{(X,Y)}^p)&\to & K_0(\M_Y^p)\\
 &\mapd{}&&\mapd{\alpha}&&\mapd{\beta}\\
 \to &  G_0(X^{p/p+1}, Y^{p/p+1})&\stackrel{\gamma}{\to}
 &K_0(\M_{(X,Y)}^{p/p+1})&\to & K_0(\M_Y^{p/p+1})
\enddiagram
we have compatible splittings $s_\alpha$ and $s_\beta$ to $\alpha$ and $\beta$, where
$s_\alpha$ is given by sending the class of an $n$-dimensional vector space over $k(x)$,
$x\in (X-Y)^p$ to the class of the free rank-$n$ module over $\O_{\bar{\{x\}}}$; $s_\beta$ is
defined similarly. This shows the first assertion.

In order to show the second assertion, it suffices to show that the map $\gamma$ above is
injective or, equivalently, that the natural map $K_1(\M_{(X,Y)}^{p/p+1})\to
K_1(\M_Y^{p/p+1})$ is surjective. This is obvious for $p=2$. Let $p=0$ and let
 \[
 a \in \bigoplus_{y \in Y^0} k(y)^\times = K_1(\M_Y^{0/1}).
 \]
Consider the semi-local ring $R \subset k(X)^\times$ of the (finitely many) generic points of
Y in X. We find a unit $\alpha \in R^\times$ reducing to $a$. By the last proposition,
$\alpha$ defines an element in $K_1(\M_{(X,Y)}^{0/1})$ which maps to $a$.

It remains the case $p=1$. Let $y$ be a codimension~$1$ point of $Y$ and $a \in k(y)^\times$.
Choose a reduced irreducible subscheme $Z$  of codimension~$1$ in $X$ which intersects $Y$
properly. Consider the semi-local ring $R\subset k(Z)$ of $Y \cap Z$ in $Z$. We find a unit
$\alpha \in R^\times$, which reduces to $a$ in $y$ and to $1$ in all other points of $Y\cap
Z$. By the last proposition, $\alpha$ defines an element of $K_1(\M_{(X,Y)}^{1/2})$ which
maps to $a$.
\end{demo}

\begin{lemma} \label{3.6}
\[
E_1^{2,-2}(X,Y) \cong \bigoplus_{x \in (X-Y)^2} \Z\; .
\]
\end{lemma}

\begin{demo}{Proof:} By proposition \ref{surj}, we have long exact sequences of the $E_1$-terms
associated to our three spectral sequences. In particular, we have an exact sequence
\[
0=E_1^{2,-3}(Y) \to E_1^{2,-2}(X,Y) \to E_1^{2,-2}(X|Y) \to E_1^{2,-2}(Y)=0,
\]
which shows the statement, since
\[
E_1^{2,-2}(X|Y)\cong K_0(\M_{(X,Y)}^{2/3})\cong \bigoplus_{x \in (X-Y)^d} \Z
\]
by proposition \ref{ausrech}.
\end{demo}

\begin{lemma} \label{3.7}
\[
E_1^{1,-2}(X,Y)= \ker \Big( \ker \big(\bigoplus_{x \in (X-Y)^1} k(x)^\times \to
\bigoplus_{y\in Y^1} \Z \big) \lang \bigoplus_{y \in Y^1} k(y)^\times\Big).
\]
\end{lemma}

\begin{demo}{Proof:}
By definition, $E_1^{1,-2}(X,Y)=G_1(X^{1/2},Y^{1/2})$. Therefore we have an exact sequence
\[
K_2(\M_Y^{1/2})\to E_1^{1,-2}(X,Y) \to K_1(\M_{(X,Y)}^{1/2}) \to K_1(\M_Y^{1/2}).
\]
All codimension~$1$ points of $Y$ have finite residue fields, and since $K_2$ of finite
fields is zero, the left hand term of the sequence vanishes. This, and the calculation of
the other terms (prop.~\ref{ausrech}) shows the statement.
\end{demo}

Lemmas \ref{3.6} and \ref{3.7} imply the

\begin{prop}\label{chow0descrklein}
The group $\CH_0(X,Y)$ is the quotient of $\ds\bigoplus_{x\in (X-Y)^2}\Z$ by the image of
\[
\ker \Big( \ker \big(\bigoplus_{x \in (X-Y)^1} k(x)^\times \to \bigoplus_{y \in Y^1} \Z \big)
\lang  \bigoplus_{y \in Y^1} k(y)^\times \Big).
\]
\end{prop}

We see that $\CH_0(X,Y)=E_1^{2,-2}(X,Y)/d( E_1^{1,-2}(X,Y))$ is the quotient of the group of
zero-cycles on $X-Y$ modulo a relation which is finer than rational equivalence. The
reciprocity homomorphism (to be constructed in section~\ref {mainsect} below) will send the
class of a closed point $P\in X-Y$ to its Frobenius automorphism in $\tilde
\pi_1^t(X,Y)^\ab$. Of course, we will have to show that this gives a well-defined
homomorphism.

\section{Finiteness of \boldmath $\CH_0(X,Y)$} \label{finsect}
The aim of this section is to detect relations in $\CH_0(X,Y)$ and to use these relations to
show that $\CH_0(X,Y)$ is finite, when $X$ is an arithmetic surface. We keep the assumption
that $X$ is a connected regular two-dimensional scheme such that all residue fields of
closed points are finite and that $Y$ is the support of a divisor on $X$. We arrange the
$E_1$-terms $E_1(X,Y)$, $E_1(X|Y)$ and $E_1(Y)$ of the spectral sequences introduced in the
last section as a (up to sign convention) double complex in the following way:

\begin{footnotesize}
\diagram{ccccccccc}
 0&&0&&0&&0&&0\\
 \mapu{}&&\mapu{}&&\mapu{}&&\mapu{}&&\mapu{}\\
 0 & \to &E_1^{2,-2}(X,Y)& \stackrel{\sim}{\to} &E_1^{2,-2}(X|Y)& \to &0&\to &0\\
 \mapu{}&&\mapu{}&&\mapu{}&&\mapu{}&&\mapu{}\\
 0 & \to& E_1^{1,-2}(X,Y)& \to &E_1^{1,-2}(X|Y)& \to &E_1^{1,-2}(Y)&\to
 & E_1^{1,-1}(X,Y)\\
\mapu{}&&\mapu{}&&\mapu{}&&\mapu{}&&\mapu{}\\
 E_1^{0,-3}(Y) & \to& E_1^{0,-2}(X,Y)& \to &E_1^{0,-2}(X|Y)&
 \twoheadrightarrow &E_1^{0,-2}(Y)&\to
 & E_1^{0,-1}(X,Y)\\
 \mapu{}&&\mapu{}&&\mapu{}&&\mapu{}&&\mapu{}\\
 0&&0&&0&&0&&0
\enddiagram
\end{footnotesize}
The rows are exact by proposition \ref{surj}. The vertical maps are the differentials of the
various spectral sequences. The term $E_1^{1,-3}(Y)$ vanishes since all residue fields of
codimension~$1$ points of $Y$ are finite, and $K_2(\hbox{finite field})=0$. All other zeros
are due to dimension reasons.

\begin{lemma} \label{surjclaim}
The map $E_1^{0,-2}(X|Y) \to E_1^{0,-2}(Y)$ is surjective.
\end{lemma}

\begin{demo}{Proof:}
By proposition~\ref{ausrech}, we know that
\[
E_1^{0,-2}(X|Y)=K_2(\M_{(X,Y)}^{0/1})= \ker \Big( K_2(k(X)) \to \bigoplus_{y \in Y^0}
k(y)^\times \Big)
\]
and
\[
E_1^{0,-2}(Y)=K_2(\M_Y^{0/1})= \bigoplus_{y \in Y^0} K_2(k(y))\;.
\]
Now recall the filtration on $K_2(F)$ of a discrete valuation field $F$ from
section~\ref{classfieldsect}.5.
For each point $y \in Y^0$ we denote by $k(X)_y$ the completion of $k(X)$ at the discrete
valuation associated with $y$. The usual approximation lemma for a finite set of discrete
valuations shows that the map
\[
K_2(k(X)) \to \bigoplus_{y \in Y^0} K_2(k(X)_y)/U^1K_2(k(X)_y)
\]
is surjective. We obtain the required statement from snake lemma applied to the commutative
diagram with exact columns

\diagram{ccc}
 0&  &0\\
 \mapd{}&&\mapd{}\\
 E_1^{0,-2}(X|Y)&\lang & \ds\bigoplus_{y\in Y^0} K_2(k(y))\\
 \mapd{}&&\mapd{}\\
 K_2(k(X))&\lang &\ds\bigoplus_{y \in Y^0} K_2(k(X)_y)/U^1K_2(k(X)_y)\\
  \mapd{}&&\mapd{}\\
 \ds \bigoplus_{y \in Y^0} k(y)^\times &\longeq& \ds\bigoplus_{y \in Y^0} k(y)^\times\\
   \mapd{}&&\mapd{}\\
   0&&0
\enddiagram
\end{demo}

\noindent
Now consider the spectral sequences associated to the considered double complex with respect
to the vertical and horizontal filtration. Since all rows are exact, the $E_1$-terms of the
first spectral sequence vanish and therefore the total complex of the double complex is
exact. The second spectral sequence degenerates at $E_4$ and its $E_1$-terms are just the
$E_2$-terms of the three spectral sequences we started with. The last lemma and the
exactness of the rows show that $E_1^{0,-2}(Y) \to E_1^{0,-1}(X,Y)$ is the zero map.
Therefore the differential $d_3$ starting from the $E_3$-term at the place $E_1^{2,-2}$
vanishes (this is a little bit confusing, but it also would not be very helpful to introduce
an explicit enumeration of the double complex). Therefore already the $E_3$-term at this
position is trivial. Explicitly this means, that $d_2$ induces  a bijection (of $E_2$-terms
of the spectral sequence associated with the double complex)

\medskip $ \ker\Big(E_2^{2,-2}(X,Y) \to E_2^{2,-2}(X|Y)\big) \liso$

\begin{flushright}
$ H\big(E_2^{1,-2}(X|Y) \to E_2^{1,-2}(Y) \to E_2^{1,-1}(X,Y)\big)$
\end{flushright}
and we get the

\begin{lemma} \label{fincrit}
$\CH_0(X,Y)$  is finite if $E_2^{2,-2}(X|Y)$ and $E_2^{1,-2}(Y)$ are finite.
\end{lemma}

Recalling the notation
\[
\SK_i(Y):=\coker \big(\bigoplus_{y \in Y^0} K_{i+1}(k(y)) \to \bigoplus_{y \in Y^1} K_i(k(y)
\,\big),
\]
for $i=0,1$, we obtain the following characterization of the group $\CH_0(X,Y)$ by a larger
set of relations than that given in proposition~\ref{chow0descrklein}.

\begin{theorem}\label{chow0descr}
$\CH_0(X,Y)$ is the quotient of $\ds\bigoplus_{x\in (X-Y)^2}\Z$ by the image of the group
\[
\ker \Big( \ker \big(\bigoplus_{x \in (X-Y)^1} k(x)^\times \to \SK_0(Y) \big) \lang \SK_1(Y)
\Big)
\]
under the divisor map.
\end{theorem}

\begin{demo}{Proof:}
By proposition~\ref{chow0descrklein}, $\CH_0(X,Y)$ is the quotient of $\ds\bigoplus_{x\in
(X-Y)^2}\Z$ by the image of
\[
\ker \Big( \ker \big(\bigoplus_{x \in (X-Y)^1} k(x)^\times \to \bigoplus_{y \in Y^1} \Z \big)
\lang  \bigoplus_{y \in Y^1} k(y)^\times \Big).
\]
Let $a$ be an element in
\[
\ker \big(\bigoplus_{x \in (X-Y)^1} k(x)^\times \to \bigoplus_{y\in Y^1} \Z \big)
=E_1^{1,-2}(X|Y).
\]
The image of $a$ in $E_1^{2,-2}(X|Y)=E_1^{2,-2}(X,Y)$ (i.e.\ the divisor of $a$) depends only
on the class of $a$ modulo the image of $E_1^{0,-2}(X|Y)$. If $a$ maps to zero in $\SK_1(Y)$,
the surjectivity of $E_1^{0,-2}(X|Y) \to E_2^{0,-2}(Y)$ implies the existence of an  element
$b\in E_1^{0,-2}(X|Y)$ with
\[
a-db \in E_1^{1,-2}(X,Y)=\ker (E_1^{1,-2}(X|Y) \lang E_1^{1,-2}(Y)).
\]
This shows that $\CH_0(X,Y)$ is the quotient of $\ds\bigoplus_{x\in (X-Y)^2}\Z$ by the image
of
\[
\ker \Big( \ker \big(\bigoplus_{x \in (X-Y)^1} k(x)^\times \to \bigoplus_{y \in Y^1} \Z \big)
\lang  \SK_1(Y) \Big).
\]
Now assume that $a$ is an element of
\[
\ker \big(\bigoplus_{x \in (X-Y)^1} k(x)^\times \to \SK_0(Y) \big).
\]
Since the map $K_2(k(X)) \to \bigoplus_{y \in Y^0} k(y)^\times$ is surjective, there exists
an element $b \in K_2(k(X))$ such that $a-db$ maps to zero already in $\bigoplus_{y \in Y^1}
\Z$, and is therefore an element in $E_1^{1,-2}(X|Y)$. Note that the divisors of $a$ and
$a-db$ (considered as elements in $\bigoplus_{x \in X^2} \Z$) coincide in all components $x
\in (X-Y)^2$. If $b' \in K_2(k(X))$ is another such element, then
\[
b-b' \in \ker \big( K_2(k(X)) \lang \bigoplus_{y \in Y^0} k(y)^\times\big)= E_1^{0,-2}(X|Y)
\]
and therefore $a -db$ and $a-db'$ (being elements in $E_1^{1,-2}(X|Y)$ have the same image
in $\SK_1(Y)=E_1^{1,-2}(Y)/d(E_1^{0,-2}(Y))$. We can therefore extend the map
\[
\ker \big(\bigoplus_{x \in (X-Y)^1} k(x)^\times \to \bigoplus_{y \in Y^1} \Z \big) \lang
\SK_1(Y)
\]
to a well-defined map
\[
\ker \big(\bigoplus_{x \in (X-Y)^1} k(x)^\times \to \SK_1(Y) \big) \lang \SK_1(Y)
\]
and the divisor of every element in the kernel of the last map gives a relation in
$\CH_0(X,Y)$.
\end{demo}

\begin{theorem} If $X$ is a connected, regular, two-dimensional scheme, flat and proper over
$\Spec(\Z)$, then the group $\CH_0(X,Y)$ is finite.
\end{theorem}

\begin{demo}{Proof:}
By lemma \ref{fincrit}, it suffices to show that $E_2^{1,-2}(Y)$ and $E_2^{2,-2}(X|Y)$ are
finite. By definition,
\[
E_2^{1,-2}(Y)=\SK_1(Y).
\]
If $Y$ would be just one regular divisor then, by Moore's theorem (\cite{Mo}, \cite{Mi}
cor.16.2), this group would be trivial if $Y$ is horizontal and isomorphic to the finite
group $\Gamma(Y,\G_m)$ if $Y$ is vertical. In the general situation we have a surjection
\[
\prod_{i=1}^r \SK_1(\tilde Y_i) \lang E_2^{1,-2}(Y)
\]
where $\tilde Y_1, \ldots, \tilde Y_r$ are the irreducible components of the normalization
of $Y$. This shows the finiteness of $E_2^{1,-2}(Y)$.

The finiteness of $E_2^{2,-2}(X|Y)$ will follow from the finiteness of
$E_2^{2,-2}(X)=\CH_0(X)$ (see \cite{K-S1}) and the injectivity of the natural map
$$E_2^{2,-2}(X|Y) \to E_2^{2,-2}(X),$$ which can be seen as follows: By
proposition~\ref{ausrech}, we have the exact sequence
 \[
0 \to E_1^{1,-2}(X|Y)\to \bigoplus_{x \in (X-Y)^1}k(x)^\times \to \bigoplus_{y \in Y^1} \Z
\to 0.
\]
Let $\sum_{i=1}^n P_i$, $P_i \in (X-Y)^2$, $i=1,\ldots,n$ represent an element of the kernel
of $E_2^{2,-2}(X|Y) \to E_2^{2,-2}(X)$. Then there exist points $x_1,\ldots, x_m \in X^1$ and
rational functions $f_j\in k(x_j)^\times$, $j=1,\ldots, m$, such that $\sum_{i=1}^n P_i =
\sum_{j=1}^m \div (f_j)$. After reordering we may assume that $x_1,\ldots, x_r\in Y^0$ and
$x_{r+1},\ldots,x_m \in (X-Y)^1$. By the same argument as in the proof of
lemma~\ref{surjclaim}, we find an element $a\in K_2(k(X))$ which has image (i.e.\ tame
symbol) $f_j^{-1}$ in $k(x_j)^\times$, $j=1,\ldots,r$ and has image $1$ in $k(y)^\times$ for
every other point $y\in Y^0$. Then
\[
\sum_{j=1}^m f_j + da \in \bigoplus_{x \in X^1}k(x)^\times
\]
has the same divisor $\sum_{i=1}^n P_i$ and nontrivial components only at points $x \in
(X-Y)^1$. We conclude that the class of $\sum_{i=1}^n P_i$ in $E_2^{2,-2}(X|Y)$ is zero. This
finishes the proof of the theorem.
\end{demo}

Furthermore, we get the
\begin{theorem} \label{horver1} Let $X$ be a connected, regular, two-dimensional scheme, flat and proper over
$\Spec(\Z)$ and let $D$ be an irreducible {\em horizontal} divisor on $X$. Then the natural
map
\[
\CH_0(X,Y \cup D) \lang \CH_0(X,Y)
\]
is an isomorphism.
\end{theorem}

\begin{demo}{Proof:}
We start with the surjectivity. For this we have to show that every element of $\CH_0(X,Y)$
can be represented by a zero-cycle $\sum_i P_i$, with $P_i \in X-(Y \cup D)$. Let $P$ be a
closed point of $D-Y$. Choose an irreducible one-dimensional subscheme $Z \subset X$
containing $P$ as a regular point and which has proper intersection with $Y\cup D$. We can
find a rational function $f$ on $Z$ which has a simple zero in $P$ and is defined and
congruent $1$ at every other point of $Z \cap (Y \cup D)$. By theorem~\ref{chow0descr}, $f$
defines a relation in $\CH_0(X,Y)$ and thus the class of $P$ coincides with the class of a
linear combination of points in $X-(Y \cup D)$.

In order to show injectivity, we consider two copies of the double complex of the beginning
of this section, namely one for $Y\cup D$ and one for $Y$, and the natural map between them.
Assume that the class of a zero-cycle $\sum P_i$, $P_i \in X-(Y\cup D)$ is zero in
$\CH_0(X,Y)$. This means that it is the image of an element $a\in E_1^{1,-2}(X,Y)$. Consider
$a$ as an element in
\[
E_1^{1,-2}(X|Y)=\ker \Big( \bigoplus_{x \in (X-Y)^1} k(x)^\times \lang \bigoplus_{y\in Y^1}
\Z \Big).
\]
Changing $a$ by an appropriate element of the image of $E_1^{0,-2}(X|Y)$, we may suppose
that the component of $a$ at the generic point of $D$ is trivial, i.e.
\[
a\in \ker \Big( \bigoplus_{x \in (X-(Y\cup D))^1} k(x)^\times \lang \bigoplus_{y\in Y^1} \Z
\Big).
\]
However, after this change we have lost the property that the image of $a$ in $E_1^{1,-2}(Y)$
is zero, but it is still contained in the image of $E_1^{0,-2}(Y)$ in $E_1^{1,-2}(Y)$. Since
all $P_i$ are in $X-(Y\cup D)$, $a$ is already an element of the subgroup
\[
E_1^{1,-2}(X|Y\cup D)=\ker \Big( \bigoplus_{x \in (X-(Y\cup D))^1} k(x)^\times \lang
\bigoplus_{y\in (Y\cup D)^1} \Z \Big).
\]
Changing $a$ by the image of an element of $E_1^{0,-2}(X,Y\cup D)$, we achieve that the
image of $a$ in
\[
E_1^{1,-2}(Y\cup D)= \bigoplus_{y \in (Y\cup D)^1} k(y)^\times
\]
has a nontrivial component only at points $y \in (D-Y)^1$. Since $D$ is horizontal, we have
$\SK_1(D)=0$, and therefore, changing $a$  by an appropriate element in $E_1^{0,-2}(X|Y\cup
D)$ (whose image in $E_1^{0,-2}(Y \cup D)$ has only a nontrivial component at the generic
point of $D$), we achieve that the image of $a$ in $E_1^{1,-2}(Y\cup D)$ is trivial, and
hence $a$ lies in the subgroup $E_1^{1,-2}(X,Y\cup D)$ of $E_1^{1,-2}(X|Y\cup D)$. All
modifications of $a$ left its divisor $\sum P_i$ invariant and thus the class of this
zero-cycle is already zero in $\CH_0(X,Y\cup D)$.
\end{demo}

\section{Tame coverings} \label{tamesect}

Coverings of a regular scheme which are tamely ramified along a normal crossing divisor have
been studied in \cite{SGA1}, \cite{G-M}. A naive extension of the (valuation theoretic)
definition of tame ramification in the normal crossing case proves to be not useful in the
general situation. For example, it would not be stable under base change (cf. \cite{S1},
Example 1.3). A definition of tameness in the general situation was given in \cite{S1} and
it was shown there that it coincides with the previous one in the normal crossing case.  The
larger freedom of having the notion of tame ramification along an arbitrary Zariski-closed
subscheme is essential in the proof of the finiteness theorem~\ref{finite}, even if $Y$ {\it
is} the support of a normal crossing divisor. However, the reader who is interested in our
main result theorem~\ref{main} only in the normal crossing case might safely skip the next
lines.

\medskip
Let us recall the definition of tameness from \cite{S1}. Let $X$ be a connected noetherian
scheme, $Y \subset X$ a closed subscheme and $U=X-Y$ the open complement. For a point $y \in
Y$ we write $X_y^{\sh}$ for $\Spec(\Cal O_{X,y}^{\sh})$, where $\sh$ means strict
henselization. By abuse of notation, we write $U_y^{\sh}$ for the base change $U \times_X
X_y^{\sh}$. The scheme $U_y^{\sh}$ is empty if $y \notin \bar U$.

\begin{defi} \label{tame2}
We say that a finite \'{e}tale covering $U' \to U$ is  tamely ramified along $Y$ if for every
point $y \in Y$ such that $U_y^{\sh}$ is nonempty, the base change
\[
U' \times_U U_y^{\sh} \lang U_y^{\sh}
\]
can be dominated by an \'{e}tale covering of the form
\[
V_1 \Tcup \cdots \Tcup V_r \lang U_{y}^{\sh},
\]
such that each $V_i$ is a connected \'{e}tale Galois covering of its image in $U_y^{\sh}$ and the
degree of\/ $V_i$ over its image is prime to the characteristic of\/ $k(y)$.
\end{defi}
If $X$ is regular and $Y$ is the support of a normal crossing divisor $D=D_1 + \cdots +D_r$
on $X$, then, by \cite{S1}, prop.1.14, this condition is equivalent to the original condition
(\cite{SGA1}, XIII.2, \cite{G-M},\S2) that the extension of function fields $k(\tilde
U)|k(U)$ is tamely ramified at the discrete valuations $v_1,\ldots,v_r$ of $k(U)=k(X)$ which
correspond to $D_1,\ldots, D_r$. For any base point $* \in U$,  the tame fundamental group
$\pi_1^t(X,Y)$ is the unique quotient of $\pi_1(U)$ which classifies finite connected \'{e}tale
coverings $\tilde{U}$ of $U$ which are tame along $Y$.  We will need the following result
from \cite{S1}.

\begin{prop}{\rm (\cite{S1}, cor.2.6)} \label{normaltame}
Let $X$ be a connected regular scheme of finite type over $\Spec(\Z)$ with function field
$K$, $U$ an open subscheme of $X$ and $Y=X-U$. Let $p$ be a prime number and let $L|K$ be a
finite Galois extension of $p$-power degree. Then the normalization\/ $U_L$ of $U$ in $L$ is
\'{e}tale over $U$ and tamely ramified along $Y$ if and only if  $L|K$ is unramified at every
discrete valuation of $K$ associated to a codimension $1$ point of $X$ which either lies on
$U$ or whose closure in $X$ contains a point of residue characteristic $p$.
\end{prop}

\begin{corol} \label{hotver2}
Assume, in addition, that $X$ is proper and flat over
$\Spec(\Z)$ and let $D$ be a horizontal prime divisor on $X$ (i.e.\ $D$ is dominant over
$\Spec(\Z)$). Then the natural map
\[
\tilde\pi_1^t(X,D \cup Y)^\ab \lang \tilde \pi_1^t(X,Y)^\ab
\]
is an isomorphism.
\end{corol}

\noindent {\bf Remark:}
 The above result is rather surprising at first glance, especially in the normal
crossing case, since tameness at a discrete valuation with residue field of characteristic
zero is a void condition. The point is that any ramification along $D$ induces additional
ramification along vertical divisors.

\begin{demo}{Proof of the corollary:}
The map in question is obviously surjective. To show injectivity, we first show that every
cyclic \'{e}tale covering $\tilde U \to U=X-(Y\cup D)$ of $p$-power degree which is tamely
ramified along $Y \cup D$ extends to an \'{e}tale covering of $V=X-Y$, i.e.\ is unramified along
$D$. Our conditions imply that $D \to \Spec(\Z)$ is proper and dominant, hence surjective.
Thus $D$ contains a point of residue characteristic $p$, and the statement follows from
proposition~\ref{normaltame}. This shows that $\pi_1^t(X,D\cup Y)^\ab \to \pi_1^t(X,Y)^\ab$
is an isomorphism. Finally, we have to deal with the real places. Let $\tilde V$ be the
normalization of $V$ in the function field of $\tilde U$. We have just shown that $\tilde V
\to V$ is \'{e}tale. By \cite{Sa}, lemma 4.9 (iii), the subset of points in $V(\R)$ which split
completely in $\tilde V$ is (norm) closed and open in $V(\R)$. Since $V(\R)$ is a real
manifold and $D\cap V(\R)$ is either empty or of real codimension $2$ in $V(\R)$, the
connected components of $U(\R)=V(\R)-D(\R)$ are in one-to-one correspondence to the
connected components of $V(\R)$. Therefore the result extends also to the modified
fundamental groups.
\end{demo}

\begin{lemma} \label{tamebyrank1}
Let $D=D_1 + \cdots +D_r$ be a sum of {\em vertical} divisors on $X$ (not necessarily with
normal crossings), $Y=\supp(D)$ and  $U=X-Y$. Let $v_1,\ldots,v_r$  be the discrete
valuations of the function field $k(U)=k(X)$ which are associated with $D_1,\ldots,D_r$.
Then a finite abelian \'{e}tale covering $\tilde U \to U$ is tamely ramified along $Y$ if and
only if the extension $k(\tilde U)|k(U)$ is tamely ramified at $v_1,\ldots,v_r$.
\end{lemma}

\begin{demo}{Proof:}
In order to show the nontrivial implication, assume that $k(\tilde U)|k(U)$ is tamely
ramified at $v_1,\ldots,v_r$. We may assume that $\tilde U \to U$ is cyclic of prime power
order, say of order $p^n$. After reordering, we may assume that $D_1,\ldots,D_s$ lie over
$p$ and that $D_{s+1},\ldots,D_r$ lie over prime numbers different to $p$. Put
$V=X-\supp(D_{s+1}+\cdots +D_r)$.  Since $v_1,\ldots, v_s$ are tamely ramified, hence
unramified in $k(\tilde U)$, the theorem on the purity of the branch locus shows that the
normalization $\tilde V$ of $V$ in $k(\tilde U)$ lies \'{e}tale over $V$. Since every point of
$Y$ with residue characteristic $p$  lies on $V$, this shows that $\tilde U \to U$ is tamely
ramified along $Y$.
\end{demo}

\section{Proof of the main theorem} \label{mainsect}
Now we are going to prove our main theorem.  We change our notation for better compatibility
with the notation of section~\ref{classfieldsect} and  use the letter $\X$ (instead of $X$)
for the scheme in question. We assume for the rest of this section that $\X$ is a
two-dimensional connected regular, proper and flat scheme over $\Spec(\Z)$. Then $X= \X
\otimes_\Z \Q$ is a smooth projective curve $\Spec(\Q)$.  Let, as in
section~\ref{classfieldsect}.4, $k$ be the algebraic closure of $\Q$ in the function field
of $X$ and let $S=\Spec(\Cal O_k)$. The structural morphism $\X \to \Spec(\Z)$ factors
through $S$ and $X$ is geometrically irreducible as a variety over $k$. Let $S_f$ be the set
of closed points of $S$ and let $S_\infty$ be the set of archimedean places of the number
field $k$. For $v \in S_f$ let $Y_v = \X \otimes_S v$ be the special fibre of $\X$ over $v$.
For $v \in S_\infty \cup S_f$ let $k_v$ be the algebraic closure of $k$ in the completion of
$k$ at $v$ and $X_v=X \times_k k_v$.

Let $\D=\D_1 + \cdots +\D_r$ be a sum of {\em vertical} divisors on $\X$ and we choose a
sufficiently small open subscheme $U\subset S$ such that $\X_U$ is disjoint to
$\Zt=\supp(\D)$. Recall the id\`{e}le group
\[
I(\X/U)= \left(\prod_{v \in P_U} \SK_1(X_v) \right) \times \left(\bigoplus_{v\in U_0}
\CH_0(Y_v) \right) ,
\]
where $P_U$ denote the set of places (including the archimedean ones) which are not in $U_0$.

\medskip
Let $v \in P_U$ be nonarchimedean and $x'$ a point in $X^1$. The set of points $x \in X_v^1$
lying over $x'$ is in one-to-one correspondence to the set of closed points of
$\overline{\{x'\}} \cap Y_v$. For a closed point $y \in Y_v^1$ we consider the
two-dimensional henselian ring
\[
A=\O_{\X,y}^h
\]
and the closed subscheme $Z=\Zt \cap \Spec(A)$, which is the support of a vertical divisor
$D=D_1 + \cdots + D_s$ on $\Spec(A)$ ($Z$ is empty, if $y \in \X -\Zt$). We define  a group
$U^t(y)$ as follows. For $y \notin \Zt$, we put
\[
U^t(y)= \ker \big(
\bigoplus_{\renewcommand{\arraystretch}{0.8}\begin{array}{c} \sst x \in X_v^1 \\
\sst x \to y
\end{array}
\renewcommand{\arraystretch}{1}}
k(x)^\times \to  \Z \big)
\]
For $y \in \Zt$, we put
\[
U^t(y)=\ker \Big( \ker \big(
\bigoplus_{\renewcommand{\arraystretch}{0.8}\begin{array}{c} \sst x \in X_v^1 \\
\sst x \to y
\end{array}
\renewcommand{\arraystretch}{1}}
k(x)^\times \to  \Z \big) \lang  k(y)^\times\Big)
\]
The group $U^t(y)$ contains the subgroup
\[
U^{1}(y)=\bigoplus_{\renewcommand{\arraystretch}{0.8}\begin{array}{c} \sst x \in X_v^1 \\
\sst x \to y
\end{array}
\renewcommand{\arraystretch}{1}} U^1(k(x)^\times).
\]

\begin{defi}
We define $U^t_\Zt C(\X/U)$ as the image of the natural composite map
\[
\Bigg(\bigoplus_{\renewcommand{\arraystretch}{0.8}\begin{array}{c} \sst v \in P_U \\
\hbox{\rm\ssz arch.} \end{array}\renewcommand{\arraystretch}{1}} \bigoplus_{x \in X_v^1}
k(x)^\times \;\; \oplus \bigoplus_{\renewcommand{\arraystretch}{0.8}
\begin{array}{c} \sst v \in P_U \\
\hbox{\rm\ssz nonarch.} \end{array}\renewcommand{\arraystretch}{1}} \bigoplus_{y \in Y_v^1}
U^t(y)\Bigg) \to I(\X/U) \to C(\X/U).
\]
\end{defi}

\begin{theorem} \label{factorprop}
Let $\X'$ be the normalization of $\X$ in a finite cyclic extension of its function field
such that $\X'_U/\X_U$ is \'{e}tale. Then $\X'/\X$ is \'{e}tale over $\X-\Zt$, tamely ramified along
$\Zt$ and every real point of $\X-\Zt$ splits completely in $\X'$ if and only if the
associated character
\[
\chi_{\X'}: C(\X/U) \to \Q/\Z
\]
annihilates the subgroup $U^t_\Zt C(\X/U)$.
\end{theorem}

\begin{demo}{Proof:}
We may assume that $\X' \to \X$ is cyclic of prime power degree, say of order $p^\alpha$.
Firstly we assume that $\X'/\X$ is \'{e}tale over $\X-\Zt$, tamely ramified along $\Zt$ and that
every real point of $\X-\Zt$ splits completely in $\X'$. The results of
section~\ref{classfieldsect}.2 show that the image of
\[
\bigoplus_{\renewcommand{\arraystretch}{0.8}\begin{array}{c} \sst v \in P_U \\
\hbox{\rm\ssz arch.} \end{array}\renewcommand{\arraystretch}{1}} \bigoplus_{x \in X_v^1}
k(x)^\times
\]
in $C(\X/U)$ is annihilated by $\chi_{\X'}$. Let $v \in P_U$ be nonarchimedean and let $y
\in Y_v$. Let $A=\O_{X,y}^h$, $Z=\Zt \cap \Spec(A)$ and let $K$ be the quotient field of
$A$. We consider the modulus $m$ on $A$ with $m(z)=1$, if $m \in Z$, and $m(z)=0$ otherwise.
$\X'$ induces a cyclic covering of $\Spec(A)$, which is \'{e}tale over $\Spec(A)-Z$ and tamely
ramified along $Z$. By corollary~\ref{localtamecrit}, the corresponding character $\chi:
C(K) \to \Q/\Z$ factors through $C_m(K)$.  Consider the defining exact sequence
\[
K_2(K)\to \bigoplus_{z \in \Spec(A)^1} K_2(K_z)/U^{m(z)}K_2(K_z) \to C_m(K) \to 0
\]
For $z \notin Z$, we have $m(z)=0$ and therefore $K_2(K_z)/U^{m(z)}K_2(K_z)\cong
k(z)^\times$. If $y\notin \Zt$, i.e.\ if $Z$ is empty, we have $C_m(K)\cong C_{(0)}(K)\cong
\Z$ and $U^t(y)$ is obviously contained in the kernel of
$\bigoplus_{\renewcommand{\arraystretch}{0.8}\begin{array}{c} \sst x \in X_v^1 \\
\sst x \to y
\end{array}
\renewcommand{\arraystretch}{1}}
k(x)^\times \lang C_m(K)$.  The commutative diagram

 \diagram{ccccc}
 C(\X/U)&\to& \pi_1^\ab(\X_U)& \to & \pi_1^t(\X,\Zt)^\ab\\
 \mapu{} &&&&\mapu{}\\
 \ds\bigoplus_{\renewcommand{\arraystretch}{0.8}\begin{array}{c} \sst x \in X_v^1 \\
 \sst x \to y
 \end{array}
 \renewcommand{\arraystretch}{1}}
 k(x)^\times &\to&C_m(K)&\to &\pi_1^t(\Spec(A),Z)^\ab
 \enddiagram
shows that it suffices to show the following

\medskip\noindent
{\it Claim:}  Assume that $y \in \Zt$. Then for every continuous homomorphism $\phi: C_m(K)
\to \Q/\Z$ with finite image of $p$-power order the induced map
\[
\bigoplus_{\renewcommand{\arraystretch}{0.8}\begin{array}{c} \sst x \in X_v^1 \\
\sst x \to y
\end{array}
\renewcommand{\arraystretch}{1}}
k(x)^\times \lang C_m(K) \lang \Q/\Z
\]
annihilates the subgroup $U^t(y)$.

\medskip\noindent
Let us prove the claim.  If $y$ has residue characteristic $p$ (i.e.\ if $v|p$), $\phi$
factors through $C_{(0)}$  and we see as above that it annihilates even the bigger group
$\ker \big(
\bigoplus_{\renewcommand{\arraystretch}{0.8}\begin{array}{c} \sst x \in X_v^1 \\
\sst x \to y
\end{array}
\renewcommand{\arraystretch}{1}}
k(x)^\times \to  \Z \big)$. It remains the case that $y \in \Zt$ has residue characteristic
different to $p$. Let $C_{m,(0)}(K)=\ker (C_m(K) \to C_{(0)}(K)=\Z)$. We have the following
diagram

 \diagram{ccc}
 \ds\ker \big(
\bigoplus_{\renewcommand{\arraystretch}{0.8}\begin{array}{c} \sst x \in X_v^1 \\
\sst x \to y
\end{array}
\renewcommand{\arraystretch}{1}}
k(x)^\times \to  \Z \big)&\to &C_{m,(0)}(K)\\
\injd{}&&\injd{}\\
\ds\bigoplus_{\renewcommand{\arraystretch}{0.8}\begin{array}{c} \sst x \in X_v^1 \\
\sst x \to y
\end{array}
\renewcommand{\arraystretch}{1}}
k(x)^\times&\to&C_m(K)
 \enddiagram
The group $C_{m,(0)}$ is a quotient of
$$
\bigoplus_{z \in Z} U^0K_2(K_z)/U^1K_2(K_z) \cong \bigoplus_{z \in Z} K_2(k(z))
$$
The continuous homomorphism $\phi$ annihilates the image of $U^iK_2(k(z))$ for some $i$.
Since $U^iK_2(k(z))/U^{i+1}K_2(k(z))$ is annihilated by $p$ for $i\geq 1$, $\phi$ is already
trivial on the image of $U^1K_2(k(z))$. Furthermore, $U^0K_2(k(z))/U^1K_2(k(z)) \cong
K_2(k(y))$ is zero, and therefore $\phi$ factors through
$K_2(k(z))/U^0K_2(k(z))=k(y)^\times$. Thus the composite
\[
\ker \big(
\bigoplus_{\renewcommand{\arraystretch}{0.8}\begin{array}{c} \sst x \in X_v^1 \\
\sst x \to y
\end{array}
\renewcommand{\arraystretch}{1}}
k(x)^\times \to  \Z \big) \to C_{m,(0)}(K) \mapr{\phi} \Q/\Z
\]
factors through the natural map

\medskip
$\ds \ker \big(
\bigoplus_{\renewcommand{\arraystretch}{0.8}\begin{array}{c} \sst x \in X_v^1 \\
\sst x \to y
\end{array}
\renewcommand{\arraystretch}{1}}
k(x)^\times \to  \Z \big) \to$
\begin{flushright}
$\ds
 \bigoplus_{z \in Z^0} K_2(k(z))/U^0K_2(k(z))= \bigoplus_{z \in
Z^0} k(y)^\times \mapr{\Sigma} k(y)^\times.\qquad$
\end{flushright}
This proves the claim and shows that $U^t_\Zt C(\X/U)$ is annihilated by $\chi_{\X'}$.

\bigskip
Now assume that $\chi_{\X'}$ annihilates $U^t_\Zt C(\X/U)$. Let $v \in P_U$ and  $x \in
X_{v,0}$. Let $y \in Y_v$ be the unique point to which $x$ specializes. By local class field
theory and by our assumption, $\chi_{\X'}$ induces a cyclic field extension of the henselian
field $k(x)$ which is tamely ramified or unramified if $y \in \Zt$ or $y \notin \Zt$,
respectively. In order to show that $\X'/\X$ is \'{e}tale over $\X -\Zt$, if suffices to show
(purity of the branch locus) that it is unramified at every prime divisor $E$ not contained
in $\Zt$.  Since $\X'_U/\X_U$ is \'{e}tale, we may assume that $E$ is vertical lying over a place
$v \in P_U$. The branch locus is closed, and so it suffices to find a closed unramified point
on $E$. This is easy as for any closed point $y\in Z$ which is a regular point of the reduced
subscheme $Y_{v, \red}$ of the fibre $Y_v$, the required property follows from
lemma~\ref{saitolemma}. It remains to show the tameness along $\Zt$. By lemma
\ref{tamebyrank1}, it suffices to show that $\X'/\X$ is tamely ramified at all generic
points of $\Zt$.  Tameness at points of residue characteristic different to $p$ is
immediate. Let $E$ be a prime divisor contained in $\Zt$ and lying in the fibre over $v \in
P_U$, $v |p$. We will show that $\X'/\X$ is unramified along $E$. Let $y$ by a closed
regular point of $E$. Every point of $X_{v,0}$ that specializes to $y$ has a residue field
which is a henselian field of characteristic zero and residue characteristic $p$. Therefore
the associated tamely ramified character $\chi_{\X'}$ on these fields is unramified. Again
using lemma~\ref{saitolemma}, we conclude that $\X'/\X$ is \'{e}tale at $y$, and therefore also
at the generic point of $E$. Hence $\X'/\X$ is tame along $\Zt$. The remaining assertion
concerning the real points follows in a straightforward manner from the results of
section~\ref{classfieldsect}.2. This concludes the proof.
\end{demo}

\begin{defi} \quad
$ C^t_\Zt(\X/U)= C(\X/U)/U^t_\Zt(\X/U).$
\end{defi}
For an abelian group $A$ we denote by  $\Div(A)$ its maximal divisible subgroup and we put
$A_\div =A/\Div(A)$.
\begin{corol}
We obtain a natural isomorphism of finite abelian groups
\[
C^t_\Zt(\X/U)_\div
 \liso \pi_1^t(\X,\Zt)^\ab.
\]
\end{corol}

\begin{demo}{Proof:}
By the last theorem, we obtain for every integer $n$ an isomorphism
\[
C^t_\Zt(\X/U)/n
 \liso \pi_1^t(\X,\Zt)^\ab/n.
\]
If $N$ annihilates the finite group $\pi_1^t(\X,\Zt)^\ab$, we have the equality
$NC^t_\Zt(\X/U)=nNC^t_\Zt(\X/U)$ for every integer $n$ and therefore
 \[
NC^t_\Zt(\X/U) \subset \Div(C^t_\Zt(\X/U)).
\]
Since $C^t_\Zt(\X/U)/N$ is finite, the last inclusion is equality. This concludes the proof.
\end{demo}

\begin{prop} \label{exreci}
The natural map
\[
\bigoplus_{x \in (\X-\Zt)^2} \Z \lang \tilde \pi_1^t(\X,\Zt)^\ab, \quad (1 \in \Z)_x \mapsto
F_x
\]
can be lifted to a map  $\bigoplus_{x \in (\X-\Zt)^2} \Z \to C^t_\Zt(\X/U)$, and this lifting
induces a surjective homomorphism
\[
\CH_0(\X,\Zt) \surjr{} C^t_\Zt(\X/U).
\]
\end{prop}

\begin{corol} \label{isocorol}
The group $C^t_\Zt(\X/U)$ is finite. In particular, we have an isomorphism
\[
C^t_\Zt(\X/U)
 \liso \pi_1^t(\X,\Zt)^\ab.
\]
\end{corol}

\begin{demo}{Proof of proposition~\ref{exreci}:} Considering the (possibly) larger subscheme
\[
\Zt' = \bigcup_{\renewcommand{\arraystretch}{0.8}\begin{array}{c} \sst v \in P_U \\
\hbox{\rm\ssz nonarch.} \end{array}\renewcommand{\arraystretch}{1}}\!\! Y_v \;\; \supset \Zt,
\]
we assume for a moment that $\Zt=\Zt'$. Then we have a natural composite map
\[
\bigoplus_{x \in (\X-\Zt)^2} \Z = \bigoplus_{v \in U_0} \bigoplus_{y \in Y_v^1} \Z  \surjr{}
\bigoplus_{v \in U_0} \SK_0(Y_v) \lang C(\X/U).
\]
Since the map $\SK_1(X) \to \prod_{v \in P_U} \SK_1(X_v)$ has a dense image, and since the
subgroup $U^t_{\Zt} C(\X/U)$ is open in $C(\X/U)$ (it contains an open subgroup of
$k(x)^\times$ for every $x \in X_v^1$, $v\in P_U$), we obtain a surjective map
\begin{equation}\label{cyclemap1}
 \bigoplus_{v \in U_0} \bigoplus_{y \in Y_v^1} \Z  \surjr{} C^t_\Zt(\X/U).
 \end{equation}
By proposition~\ref{chow0descrklein}, we know that $\CH_0(\X,\Zt)$ is the quotient of
$\ds\bigoplus_{x\in (\X-\Zt)^2}\Z$ by the image of
\[
R=\ker \Big( \ker \big(\bigoplus_{x \in (\X-\Zt)^1} k(x)^\times \to \bigoplus_{y \in \Zt^1}
\Z \big) \lang  \bigoplus_{y \in \Zt^1} k(y)^\times \Big).
\]
In order to show that the map \eqref{cyclemap1} factors through $\CH_0(\X,\Zt)$, we have to
show that it annihilates the image of $R$. Consider the diagonal map
\[
\bigoplus_{x \in (\X-\Zt)^1} k(x)^\times \to
\bigoplus_{\renewcommand{\arraystretch}{0.8}\begin{array}{c} \sst v \in P_U \\
\hbox{\rm\ssz nonarch.} \end{array}\renewcommand{\arraystretch}{1}}\bigoplus_{x \in X_v^1}
k(x)^\times=
\bigoplus_{\renewcommand{\arraystretch}{0.8}\begin{array}{c} \sst v \in P_U \\
\hbox{\rm\ssz nonarch.} \end{array}\renewcommand{\arraystretch}{1}} \bigoplus_{y \in Y_v^1}
\bigoplus_{\renewcommand{\arraystretch}{0.8}\begin{array}{c} \sst x \in X_v^1 \\
\sst x \to y
\end{array}
\renewcommand{\arraystretch}{1}}
k(x)^\times
\]
By our temporary assumption $\Zt=\Zt'$, we can read off the vanishing of $R$ in
$C^t_\Zt(\X/U)$ directly from the definition of $U^t_\Zt C(\X/U)$.  Since $\CH_0(\X,\Zt)$ is
finite, also $C^t_\Zt(\X/U)$ is finite in this case.

\medskip Returning to the general case, assume that $\Zt \subsetneq \Zt'$. The finiteness of
$C^t_{\Zt'}(\X/U)$ and the surjection $C^t_{\Zt'}(\X/U)\to C^t_\Zt(\X/U)$ show that also
$C^t_\Zt(\X/U)$ is finite. Let $v \in P_U$ and let $y\in Y_v^1$ be a closed point not
contained in $\Zt$. Then the map
\[
\bigoplus_{\renewcommand{\arraystretch}{0.8}\begin{array}{c} \sst x \in X_v^1 \\
\sst x \to y
\end{array}
\renewcommand{\arraystretch}{1}}
k(x)^\times   \lang C^t_\Zt(\X/U)
\]
annihilates $U^t(y)$ and therefore induces a well-defined map
\[
\Z \lang C^t_\Zt(\X/U) \liso \pi_1^t(\X,\Zt)^\ab.
\]
The image of $1 \in\Z$ in $\pi_1^t(\X,\Zt)^\ab$ is just the Frobenius $F_y$. Combining these
maps with the maps for $x \in (\X-\Zt')^2$ constructed above, we obtain a surjective map
\[
\bigoplus_{x \in (\X-\Zt)^2} \Z \lang C^t_\Zt(\X/U).
\]
Trying to proceed in the same way as above in the case that $\Zt=\Zt'$, we now run into
trouble with the vertical components in $(\X -\Zt)^1$ lying over places $v \in P_U$. Let
$v\in P_U$ and let $y \in Y_v \cap \Zt$ be a point, such that (setting $A=\O_{\X,y}$ and all
other notation as above) $\Spec(A)$ has a vertical prime divisor not contained in $\Zt$. Let
$\phi: C^t_\Zt(\X/U) \to \Q/\Z$ be any homomorphism. Via the commutative diagram
 \diagram{ccccc}
 \ds\bigoplus_{x \in (\Spec(A)-Z)^1} k(x)^\times &\to C^t_\Zt(\X/U) &\liso& \pi_1^t(\X,\Zt)^\ab \\
 \eqd&&&\mapu{}\\
 \ds\bigoplus_{x \in (\Spec(A)-Z)^1} k(x)^\times &\to C_m(K)&\to& \pi_1^t(\Spec(A),Z)^\ab
 \enddiagram
we get a homomorphism with finite image $\phi_y: C_m(K) \to \Q/\Z$. Exactly the same proof as
that of the claim in the proof of theorem~\ref{factorprop} shows that $\phi_y$ and hence
also $\phi$ annihilates the image of
\[
\ker \Big( \ker \big(\bigoplus_{x \in (\Spec(A)-Z)^1} k(x)^\times \to  \Z \big) \lang
k(y)^\times \Big).
\]
The finiteness of $C^t_\Zt(\X/U)$ then shows that these elements are zero in $C^t_\Zt(\X/U)$
and we are done.
\end{demo}

\begin{prop} \label{injprop}
The homomorphism $\CH_0(\X,\Zt) \to C^t_\Zt (\X/U)$ is injective.
\end{prop}

\begin{demo}{Proof:} Let $\sum n_i P_i \in \bigoplus_{x \in (\X-\Zt)^1} \Z$ be a zero-cycle on
$\X-\Zt$ which represents an element in $\ker(\CH_0(\X,\Zt) \to C^t_\Zt (\X/U))$. By
definition, there exist elements
\[
a \in \bigoplus_{x \in X^1} k(x)^\times, \; b \in \bigoplus_{v \in U_0} k(Y_v)^\times; \; c
\in \bigoplus_{\renewcommand{\arraystretch}{0.8}\begin{array}{c} \sst v \in P_U \\
\hbox{\rm\ssz nonarch.} \end{array}\renewcommand{\arraystretch}{1}} K_2(k(X_v)),
\]
such that $\sum n_i P_i = \div(a) + \div(b)$ and, denoting the diagonal image of $a$ in
\[
\bigoplus_{\renewcommand{\arraystretch}{0.8}\begin{array}{c} \sst v \in P_U \\
\hbox{\rm\ssz nonarch.} \end{array}\renewcommand{\arraystretch}{1}}\bigoplus_{x \in X_v^1}
k(x)^\times=
\bigoplus_{\renewcommand{\arraystretch}{0.8}\begin{array}{c} \sst v \in P_U \\
\hbox{\rm\ssz nonarch.} \end{array}\renewcommand{\arraystretch}{1}} \bigoplus_{y \in Y_v^1}
\bigoplus_{\renewcommand{\arraystretch}{0.8}\begin{array}{c} \sst x \in X_v^1 \\
\sst x \to y
\end{array}
\renewcommand{\arraystretch}{1}}
k(x)^\times
\]
by the same letter, we have
\[
a -\delta c \in \bigoplus_{\renewcommand{\arraystretch}{0.8}\begin{array}{c} \sst v \in P_U \\
\hbox{\rm\ssz nonarch.} \end{array}\renewcommand{\arraystretch}{1}} \bigoplus_{y \in Y_v^1}
U^t(y).
\]
If already $a$ would lie in this subgroup, then by proposition~\ref{chow0descrklein} and by
the definition of $U^t(y)$, $\div(a)+\div(b)$ would be a relation in $\CH_0(\X,\Zt)$ and we
were done. We will achieve this by modifying $a$ by the tame symbol of a suitable element in
$K_2(k(X))$ (which does not change $\div(a)$). Consider the commutative diagram
 \diagram{ccc}
 \ds\bigoplus_{v \in P_U} K_2(k(X_v))& \lang &
\ds \bigoplus_{v \in P_U} \bigoplus_{z \in Y_v^0} K_2(k(X_v)_z)\\
 \mapu{}&&\mapu{}\\
 K_2(k(X)) &\lang&\ds \bigoplus_{\renewcommand{\arraystretch}{0.8}\begin{array}{c} \sst x \in \X^1 \\
\sst z \to P_U \end{array}\renewcommand{\arraystretch}{1}} K_2(k(X)_z)
 \enddiagram
Let $v \in P_U$,  $z \in Y_v^0$ and let $y$ be a closed point in ${\overline{\{z\}}}$. Let
$K$ be the quotient field of $A=\O_{\X,y}^h$. We consider $k(X)$ as a discrete valuation
field with the valuation associated with $z$. Then  $K_2(k(X))/U^0K_2(k(X))\cong
k(z)^\times$ and $U^0K_2(k(X))/U^1K_2(k(X))\cong K_2(k(z))$ and we obtain isomorphisms
\[
K_2(k(X))/U^1 \liso K_2(k(X)_z)/U^1 \liso K_2(k(X_v)_z)/U^1  \liso K_2(K_z)/U^1.
\]
If $z_1,\ldots,z_n$ are finitely many vertical points, the usual approximation lemma shows
that the natural map
\[
K_2(k(X)) \lang \bigoplus_{i=1}^n K_2(k(X))/U^1_{z_i} K_2(k(X))
\]
is surjective. Therefore we find an element $e \in K_2(X)$ having the same image as $c$ in
\[
\bigoplus_{v \in P_U} \bigoplus_{z \in Y_v^0} K_2(k(X_v)_z) /U^1 K_2(k(X_v)_z)
\]
Then $\delta(c-e)$ has trivial vertical components (we are only interested in those lying
over a nonarchimedean point of $P_U$) and is therefore contained in the subgroup
\[
\ker \big(
\bigoplus_{\renewcommand{\arraystretch}{0.8}\begin{array}{c} \sst x \in X_v^1 \\
\sst x \to y
\end{array}
\renewcommand{\arraystretch}{1}}
k(x)^\times \to  \Z \big)
\]
and we claim that it is contained in $U^t(y)$ for every $y\in Y_v$, $v \in P_U$. Since the
statement becomes only stronger in this way, we may temporarily  replace $\Zt$ by the
(possible) larger subscheme \quad $
\bigcup_{\renewcommand{\arraystretch}{0.8}\begin{array}{c} \sst v \in P_U \\
\hbox{\rm\ssz nonarch.} \end{array}\renewcommand{\arraystretch}{1}}\!\! Y_v \;\; \supset
\Zt.$ Then the group
\[
U^t(y)=\ker \Big( \ker \big(
\bigoplus_{\renewcommand{\arraystretch}{0.8}\begin{array}{c} \sst x \in X_v^1 \\
\sst x \to y
\end{array}
\renewcommand{\arraystretch}{1}}
k(x)^\times \to  \Z \big) \lang  k(y)^\times\Big)
\]
is nothing else but the term $E_1^{1,-2}(X,Y)$ in the double complex of the beginning of
section~\ref{finsect} with $X=\Spec(A)$ and $Y=Z =\Spec(A)\cap\Zt$. The image of  $c-e$ in
$K_2(K)$ is in the kernel of  $K_2(K) \to \bigoplus_{z\in Z} K_2(K_z)/U^1K_2(K_z)$ and
therefore, in the notation of section~\ref{finsect}, it is an element in the group
$\ker(E_1^{0,-2}(X|Y))\to E_1^{0,-2}(Y)$. Hence the image of $\delta(c-e)$ in
$E_1^{0,-2}(X|Y)$ maps to zero in $E_1^{1,-2}(Y)$ and is in $E_1^{1,-2}(X,Y)$. This shows
that $\delta(c-e)$ is in $U^1(y)$. Setting $a'=a-\delta e$, we see that $a'=a-\delta c +
\delta(c-e)$ is in $\bigoplus_{\renewcommand{\arraystretch}{0.8}
\begin{array}{c} \sst v \in P_U \\
\hbox{\rm\ssz nonarch.} \end{array}\renewcommand{\arraystretch}{1}} \bigoplus_{y \in Y_v^1}
U^t(y)$.
\end{demo}

Summing up, we have proven theorem~\ref{main}. Indeed, if $Y$ is the union of vertical
divisors, then the assertion follows from proposition~\ref{injprop} and
corollary~\ref{isocorol}. The general case follows from corollary \ref{hotver2} and
theorem~\ref{horver1}.

\section{Variants}
Let, as before, $X$ be a two-dimensional regular, connected flat and proper scheme of finite
type over $\Spec(\Z)$ and let $Y$ be an {\em arbitrary} proper closed subscheme in $X$. In
view of proposition~\ref{chow0descrklein}, it seems reasonable to make the following

\begin{defi}
 The group $\CH_0(X,Y)$ is the quotient of $\ds\bigoplus_{x\in
(X-Y)_0}\Z$ by the image of
\[
\ker \Big( \ker \big(\bigoplus_{x \in (X-Y)_1} k(x)^\times \to \bigoplus_{y \in Y_0} \Z \big)
\lang  \bigoplus_{y \in Y_0} k(y)^\times \Big).
\]
\end{defi}

Then we have the following generalization of theorem~\ref{main}.

\begin{theorem}
There exists a natural isomorphism of finite abelian groups
\[
\rec_{X,Y}: \CH_0(X,Y) \lang \tilde \pi_1^t(X,Y)^\ab\;.
\]
\end{theorem}

\begin{demo}{Proof:} Let $Y=D \tcup E$, where $D$ is a union of divisors and $E$ is union of
closed points. Let $\tilde X$ be the blow-up of $X$ in $E$ and let $\tilde Y$ be the
pre-image of $Y$ in $\tilde X$. By theorem~\ref{finite}, the natural homomorphism
\[
\tilde\pi_1^t(\tilde X,\tilde Y) \lang \tilde\pi_1^t(X,Y)
\]
is an isomorphism since both groups only depend on $U=\tilde X-\tilde Y=X-Y$. In order to
prove the theorem it therefore suffices to show that the natural isomorphism
\[
\bigoplus_{x \in \tilde X -\tilde Y} \Z \liso \bigoplus_{x \in X - Y} \Z
\]
induces an isomorphism $\CH_0(\tilde X,\tilde Y) \liso \CH_0(X,Y)$. Let $a \in \bigoplus_{x
\in U_1} k(x)^\times$ be a relation in $\CH_0(\tilde X,\tilde Y)$. Then it is also relation
in $\CH_0(X,Y)$. On the other hand, let $Z$ be the pre-image of a point in $E$. Then $Z$ is a
projective line over a finite field, hence $\CH_0(Z) \cong \Z$. If $a$ is a relation in
$\CH_0(X,Y)$ it therefore follows from theorem~\ref{chow0descr} that it also defines a
relation in $\CH_0(\tilde X,\tilde Y)$. This concludes the proof.
\end{demo}

Without proof, we give yet another variant of theorem~\ref{main} which describes the full
group $\pi_1^t(X,Y)^\ab$ and not its quotient $\tilde\pi_1^t(X,Y)^\ab$. Of course, this
discussion is void, if $X$ has no real points. We define a modified relative Chow group by
allowing only such relations $\div (f)$ of functions $f$ which are positive at every
real-valued point. More precisely, for a point $x \in X_1$ we put
\[
k(x)_+^\times=\{f \in k(x)^\times\,|\, \iota(f)>0 \hbox{ for every embedding } \iota: k(x)
\to \R \}.
\]
If the global field $k(x)$ is of positive characteristic or a totally imaginary number
field, then $k(x)^\times_+ =k(x)^\times$.

\begin{defi}
The group $\widetilde \CH_0(X,Y)$ is the quotient of $\ds\bigoplus_{x\in (X-Y)_0}\Z$ by the
image of
\[
\ker \Big( \ker \big(\bigoplus_{x \in (X-Y)_1} k(x)^\times_+ \to \bigoplus_{y \in Y_0} \Z
\big) \lang  \bigoplus_{y \in Y_0} k(y)^\times \Big).
\]
\end{defi}

We then have the

\begin{theorem}
There exists a natural isomorphism of finite abelian groups
\[
\widetilde \rec_{X,Y}: \widetilde \CH_0(X,Y) \lang \pi_1^t(X,Y)^\ab\;.
\]
\end{theorem}

\vfill \noi{{\sc Alexander Schmidt}, Mathematisches Institut, Universit\"{a}t Heidelberg, Im
Neu\-enheimer Feld 288, 69120 Heidelberg, Deutschland\\ e-mail:
schmidt@mathi.uni-heidelberg.de}
\end{document}